\documentclass{conm-p-l}
\usepackage{amssymb}
\input xy
\xyoption{all}
\newtheorem{theorem}{Theorem}[section]
\newtheorem{theorem*}{Theorem}
\newtheorem{lemma}[theorem]{Lemma}
\newtheorem{proposition}[theorem]{Proposition}
\newtheorem{cor}[theorem]{Corollary}
\newtheorem{cor*}{Corollary}
\newtheorem{conjecture}{Conjecture}

\theoremstyle{definition}
\newtheorem{definition}[theorem]{Definition}
\newtheorem{definition*}{Definition}
\newtheorem{example}[theorem]{Example}

\theoremstyle{remark}
\newtheorem{remark}[theorem]{Remark}
\newtheorem{remarks}[theorem]{Remarks}
\newtheorem{remarks*}{Remarks}

\numberwithin{equation}{section}
\newcommand{\K}{\mathcal K}

\newcommand{\Ch}{{\operatorname{Ch}}}

\newcommand{\SO}{\operatorname{SO}}

\newcommand{\rank}{\operatorname{rank}}
\newcommand{\tr}{\operatorname{tr}}

\newcommand{\Todd}{{\operatorname{Todd}}}
\newcommand{\nc}{\newcommand}
\nc{\hM}{\widehat M}
\nc{\cC}{\mathcal C}
\nc{\cL}{\mathcal L}
\nc{\cA}{\mathcal A}
\nc{\bbC}{\mathbb C}
\nc{\bbR}{\mathbb R}
\nc{\bbZ}{\mathbb Z}
\nc{\bbN}{\mathbb N}
\nc{\bbQ}{\mathbb Q}
\nc{\bbH}{\mathbb H}
\nc{\Da}{\Delta}
\nc{\da}{\delta}
\nc{\ta}{\theta}
\nc{\za}{\zeta}
\nc{\A}{\mathcal A}
\nc{\Om}{\Omega}
\nc{\Hj}{{\bar H}_{(2)}^j(\hM)}
\nc{\Hdot}{{\bar H}_{(2)}^{\bullet}(\hM)}
\nc{\Oj}{\Om_{(2)}^j(\hM)}
\nc{\vp}{\varphi}
\nc{\tildvp}{\tilde\vp}
\nc{\hN}{\widehat N}
\nc{\Cj}{C^j(\hK)}
\nc{\HjK}{{\bar H}_{(2)}^j(\hK)}
\nc{\HdotK}{{\bar H}_{(2)}^{\bullet}(\hK)}
\nc{\al}{\alpha}
\nc{\M}{{\mathbb H}}
\nc{\Tau}{T}
\def\Ga{\Gamma}

\def\al{\alpha}

\def\Om{\Omega}

\def\rank{\operatorname{rank}\,}

\def\sup{\operatorname{sup}}
\def\inf{\operatorname{inf}}

\def\C{\mathbb C}
\def\R{\mathbb R}
\def\Z{\mathbb Z}

\def\N{\mathbb N}

\def\index{\operatorname{index}}

\def\A{{\mathcal A}}
\def\B{{\mathcal B}}

\def\Cl{{\mathcal Cl}}
\def\F{{\mathcal F}}

\def\tr{\operatorname{tr}}

\def\np{\operatorname{\not\!\partial}}

\def\U{{\mathcal U}}
\def\<{\langle}
\def\>{\rangle}

\def\ldg{\ell^2\Gamma}
\def\Uu{{\rm\bf U}(1)}
\nc{\uE}{{\underline E}}
\begin{document}

\title[Heat kernels and completions...]{Heat kernels and the range of the trace on completions
of twisted group algebras}
\author{Varghese Mathai}
\address{Department of Mathematics, University of Adelaide, Adelaide  
5005, Australia}
\date{\today}
\email{vmathai@maths.adelaide.edu.au}
\begin{abstract}
Heat kernels are used in this paper to express the analytic index of projectively
invariant Dirac type operators on $\Gamma$ covering spaces of compact
manifolds, as elements in the K-theory of certain unconditional completions of the
twisted group algebra of $\Gamma$.
This is combined with V. Lafforgue's results in the untwisted case,  to
compute the range of the trace on the K-theory of these algebras, under
the hypothesis that $\Gamma$ is in the class C' (defined by V. Lafforgue).
\end{abstract}
\maketitle
%------------------------------------------------------------
%------------------------------------------------------------
\section*{Introduction}
%------------------------------------------------------------
For $\Gamma$ a torsion-free discrete group, one formulation of 
a standing conjecture of Kaplansky and Kadison states that the range of the canonical trace
on the $K$-theory of the reduced $C^*$-algebra of $\Gamma$, is contained in the integers.
However, if we twist the convolution by a multiplier (i.e. a normalized $\Uu$-valued 2-cocycle on $\Gamma$) 
then this is no longer true, as shown
by Pimsner-Voiculescu \cite{PV} and Rieffel \cite{Rieff}, who computed the precise range of the canonical trace
on the twisted group $C^*$-algebra of ${\bbZ^2}$ (which turns out to be the noncommutative torus).
The author with Marcolli in \cite{MM} settled the case of surface groups by identifying the range of the canonical trace on the twisted group $C^*$-algebra.
In the present paper, we study the range of the trace on the $K$-theory of  {\em good} unconditional completions  $\cA(\Gamma, \sigma)$ of the twisted group algebra $\bbC(\Gamma,\sigma)$ (see section \ref{RD}) - an example of such a completion is the $\ell^1$ completion of $\bbC(\Gamma, \sigma)$.
Our approach is to study a twisted analogue of the assembly map, viewed as a homomorphism 
%%% 
\begin{equation}\label{twistedassembly} 
\mu_{\sigma}^{\cA} : K^\Gamma_* (\underline{E}\Gamma) \rightarrow 
K_*(\cA(\Gamma, \sigma)),
\end{equation} 
%%% 
whenever the Dixmier-Douady invariant $ \delta(\sigma)$ of the multiplier $\sigma$ on $\Gamma$ is trivial. 
The map $\mu_\sigma^\cA$ is a twisted version of the assembly map defined by Lafforgue, 
and the definition uses Lafforgue's Banach $KK$-theory. We use the method of heat kernels to study 
an analytic version of the map \eqref{twistedassembly}, called the twisted analytic 
Baum-Connes map in the paper, and a standard index theorem in section \ref{equivalence} establishes that both 
definitions are equivalent. 
Using fundamental results in \cite{Laff}, we prove in Theorem \ref{main1} that if
$\Gamma$ is a discrete group  in Lafforgue's class $\cC'$, then the twisted assembly map $\mu_\sigma^\cA$ is an isomorphism.
The class $\cC'$ will be described later, but we mention here that it contains all discrete subgroups of connected Lie groups, word hyperbolic groups and amenable groups. Together with a twisted version of an $L^2$-index theorem given in \cite{Ma2} it is then a straightforward corollary to obtain a formula for the range of the trace on the $K$-theory of  $\cA(\Gamma, \sigma)$
in terms of classical characteristic classes on the the classifying space $\underline{B}\Gamma$
for proper actions, as explained in section \ref{sect:range}. 
Although this formula is computationally challenging, it can be explicitly computed 
 in low dimensional cases, e.g. in the case 
when $\Gamma$ is torsion-free and $B\Gamma$ is a smooth compact oriented 
manifold of dimension less than or equal
to 4, which is done in section \ref{3&4} of the paper. 
This generalizes earlier results of  \cite{CHMM}, \cite{MM} ,  \cite{PV}, \cite{Rieff}, \cite{BC}.  

If in addition $\Gamma$ has the Rapid Decay property (property RD), then we can choose the good unconditional completion 
$\cA(\Gamma,\sigma)$ to be the Sobolev completion $H^s(\Gamma, \sigma)$ of $\bbC(\Gamma, \sigma)$, for $s$ large and $H^s(\Gamma, \sigma)$ being a dense $*$-subalgebra of the reduced $C^*$-algebra $C^*_r(\Gamma, \sigma)$, stable under the holomorphic functional calculus, the twisted assembly map reduces to the usual twisted assembly map cf. \cite{Ma},
\begin{equation}\label{twistedassembly2} 
\mu_{\sigma}: K^\Gamma_i (\underline{E}\Gamma) \rightarrow 
K_i(C^*_r(\Gamma, \sigma)). 
\end{equation} 
We prove in Theorem \ref{main1} that if $\Gamma$ is a discrete group in the class 
$\cC'$ and $\Gamma$ has property RD, then the twisted 
assembly map $\mu_\sigma$ in (\ref{twistedassembly2}) is an isomorphism.
The groups that have property RD include all finitely generated 
groups of polynomial growth, all finitely generated free groups, 
all word hyperbolic groups, and certain property $(T)$ groups such as 
cocompact lattices in ${\bf SL}(3, \mathbb F)$ or the exceptional 
group ${\bf E}_6$, where $\mathbb F$ is a non-discrete locally compact field or the quaternions. All of these groups
are also in the class $\cC'$. Using the earlier mentioned procedure, we also  obtain a formula for the range of the trace on the $K$-theory of  $C^*_r(\Gamma, \sigma)$
in terms of classical 
characteristic classes on $\underline{B}\Gamma$.

The last section of the paper is devoted to studying the degree one of the assembly map, and following the construction of Natsume \cite{Natsume} and Valette et. al.  in \cite{Bettaieb-Matthey}, we determine the explicit generators of $K_1(\cA(\Gamma, \sigma))$ and $K_1(C^*_r(\Gamma, \sigma))$, whenever $\Gamma$ is a torsion-free cocompact Fuchsian group. 

The appendix, written by Indira Chatterji, establishes useful results 
on the twisted rapid decay property for $(\Gamma, \sigma)$ that are used in the 
text. One interesting result there is that if $\Gamma$ has property RD, 
then it has the twisted RD property for any multiplier $\sigma$ on $\Gamma$.
This means in particular that we {\em do not} have to appeal to the Baum-Connes 
conjecture {\em with coefficients}, which is a technical improvement of results in \cite{Ma}.
The author thanks Indira Chatterji for helpful discussions.
 
%%%%%%%%%%%%%%%%%%%%%%%%%%%%%%%%%%%%%%%%%%%%%%%%%%%%%%%%%%%%%%.

%%%%%%%%%%%%%%%%%%%%%%%%%%%%%%%%%%%%%%
\section{Basics}\label{basics}
%%%%%%%%%%%%%%%%%%%%%%%%%%%%%%%%%%%%%%%%%%
In this section $\sigma$ is a multiplier on $\Gamma$ a discrete group,  
that is
a map $\sigma:\Gamma\times\Gamma\to\Uu$ satisfying the
following identity for all $\gamma, \mu,\delta\in \Gamma$:
\begin{itemize}
\item[(1)]$\sigma(\gamma,\mu)\sigma(\gamma\mu,\delta)=\sigma(\gamma,\mu\delta)\sigma(\mu,\delta)$.
\item[(2)]$\sigma(\gamma,1)=\sigma(1,\gamma)=1$.
\end{itemize}
Recall that the \emph{Dixmier-Douady invariant} of a  
multiplier $\sigma$ is the cohomology class $\delta(\sigma)\in  
H^3(\Gamma,\bbZ)$, the image of $[\sigma]$ obtained under the map  
$\delta$ arising in the long exact sequence in cohomology derived from  
the short exact sequence of coefficients
$$0\to\bbZ\to\bbR\to\Uu\to 0.$$
We denote by $E\Ga$ the universal cover of $B\Ga$, the classifying space for $\Ga$. The following lemma will be used later.
\begin{lemma}\label{twistedIdentities} Let $\alpha\in Z^2(B\Ga,\bbR)$ and $X\subset E\Ga$ a cocompact $\Ga$-space. Then there is a map $\varphi:\Ga\to C_0(X)$ such that;
\begin{itemize}
\item[(i)]$\varphi_{\gamma}(x)+\varphi_{\mu}(\gamma x)-\varphi_{\mu\gamma}(x)$ is independent of $x\in X$.
\item[(ii)]There is $x_0\in X$ such that $\varphi_{\gamma}(x_0)=0$ for any $\gamma\in\Ga$.
\item[(iii)] $\lambda(\gamma,\mu)= \varphi_{\gamma}(\mu x_0)$  is an $\bbR$-valued 2 cocycle that is cohomologous 
to $\alpha$.
%$\sigma'(\gamma,\mu)=e^{-i\varphi_{\gamma}(\mu x)}$ is a multiplier 
%which is cohomologous to the multiplier $\sigma = e^{i\alpha}$.
\end{itemize}
\end{lemma}
\begin{proof}
%Since $\sigma$ has trivial Dixmier-Douady invariant, there is an $\alpha\in Z^2(B\Ga,\bbR)$ 
%such that $[e^{\alpha}]=[\sigma]$ (where by abuse of notation we write $[\sigma]\in H^2(B\Ga,\Uu)$ 
%for the image of $[\sigma]\in H^2(\Ga,\Uu)$ under the canonical isomorphism $H^2(\Ga,\Uu)
%\simeq H^2(B\Ga,\Uu)$, see \cite{Brown}). 
Let $p:E\Ga\to B\Ga$ be the canonical projection and take a lift $\tilde{\alpha}=p^*(\alpha)\in Z^2(E\Ga,\bbR)$. Since $E\Ga$ is contractible, there is a $\Lambda\in C^1(E\Ga,\bbR)$ such that $\tilde{\alpha}=d\Lambda$. By definition of $\tilde{\alpha}$, we have
$$0=\gamma^*\tilde{\alpha}-\tilde{\alpha}=d(\gamma^*\Lambda-\Lambda)\ \ \hbox{ for any }\gamma\in\Ga.$$
The element $\eta_{\gamma}=\gamma^*\Lambda-\Lambda$ hence belongs to $Z^1(E\Ga,\bbR)$, so that there exists $c_{\gamma}\in C^0(E\Ga,\bbR)$ with $\eta_{\gamma}=dc_{\gamma}$. Let us show that $\mu^*c_{\gamma}+c_{\mu}-c_{\gamma\mu}\in C^0(E\Ga,\bbR)$ is a constant. To do so, it is enough to see that $d(\mu^*c_{\gamma}+c_{\mu}-c_{\gamma\mu})=0$. We compute:
\begin{eqnarray*}dc_{\gamma\mu}&=&\eta_{\gamma\mu}=(\gamma\mu)^*\Lambda-\Lambda=\mu^*\gamma^*\Lambda-\gamma^*\Lambda+\gamma^*\Lambda-\Lambda\\
&=&\mu^*\eta_{\gamma}+\eta_{\mu}=d(\mu^*c_{\gamma}+c_{\mu})\end{eqnarray*}
Let $x_0\in X$, we now define
\begin{eqnarray*}\varphi:\Ga &\to & C_0(X)\\
\gamma &\mapsto &\varphi_{\gamma},\end{eqnarray*}
where $\varphi_\gamma(x) =  c_{\gamma}(x)-c_{\gamma}(x_0)$.
Then $\varphi$ satisfies (i) and (ii). In particular, from the $\bbR$-valued closed 2-form $\alpha$ on $B\Gamma$,
 we have produced an $\bbR$-valued group 2-cocycle $\lambda(\gamma,\mu)= \varphi_{\gamma}(\mu x_0)$ on $\Gamma$. 
 Also, the group extension corresponding to $\lambda$ can be described as follows.
Let $\Gamma^\lambda = \Gamma \times \bbR$ with product given by 
$(\gamma, r)(\gamma', r') = (\gamma\gamma', \lambda(\gamma, \gamma') +  r  + r')$. 
If $g_{ij}$ are transition functions for the principal bundle $p: E\Gamma \to B\Gamma$, 
define the lift  $\hat g_{ij} = (g_{ij}, 0) \in \Gamma^\lambda$
% in the case when $\lambda(g, g^{-1}) = 1$, and 
%more generally $\hat g_{ij} = (g_{ij}, \xi_i \xi_j)$ if this is not the case. 
Then $ t_{ijk} = \hat g_{ij} \hat g_{jk} \hat g_{ki} = \lambda( g_{ij}, g_{jk}) + \lambda(g_{ik}, g_{ki}) $
is the $\bbR$-valued Cech 2-cocycle on $B\Gamma$ that is associated to 
the $\bbR$-valued group 2-cocycle $\lambda$ on $\Gamma$.

If $\alpha|_{U_i} = d \theta_i$, 
$(\theta_i - \theta_j)|_{U_i\cap U_j} = d f_{ij}$,  $ (f_{ij} + f_{jk}
+ f_{ki}) = t_{ijk} \in \mathbb R$, then the Cech cohomology 2-cocycle corresponding to 
the de Rham closed 2-form $\alpha$
is $t$, by the well known Cech-de Rham isomorphism.

This shows that $[\alpha] = [t] = [\lambda] \in H^2(B\Gamma, \bbR)$.

\end{proof}
In the notation of Lemma \ref{twistedIdentities} above, one verifies that 
$\sigma(\gamma, \mu) = \exp(-i\varphi_\gamma(\mu x_0))$
defines a multiplier on $\Gamma$.  The map $\varphi$ is called a \emph{phase} associated to $\sigma$.

\begin{definition}
Let $A$ be a $\Gamma$-$C^*$-algebra, we denote by $\C(\Gamma, A, \sigma)$ the 
$*$-algebra of finitely supported maps from $\Gamma$ to $A$ endowed with a 
$\sigma$-twisted convolution given as follows: for all  $a, b \in A$ and $\gamma. \mu \in \Gamma$,
$$aT_{\gamma}*_\sigma bT_{\mu}=a\alpha_{\gamma}(b)\sigma(\gamma,\mu)T_{\gamma\mu},$$
where $\alpha$ denotes the action of $\Gamma$ on $A$. Here we think of elements of 
$\C(\Gamma, A, \sigma)$ as finite sums $\sum a_\gamma T_\gamma$, where $a_\gamma \in A$, 
$T_\gamma$ is a formal letter satisfying $T_{\gamma}T_{\mu} = \sigma(\gamma,\mu)T_{\gamma\mu}$, $T_\gamma a T_\gamma^* = \alpha_\gamma (a)$ and $T_\gamma^* = \sigma(\gamma, \gamma^{-1}) T_{\gamma^{-1}}$.

Given a Banach norm $\|\ \|_{B}$ on $\C(\Gamma, A, \sigma)$, we denote by $B(\Gamma,A,\sigma)$ the completion of $\C(\Gamma, A, \sigma)$ with respect to the norm $\|\ \|_{B}$.\end{definition}
In case where $A=\C$ (with a trivial $\Ga$-action) we simply write $\bbC(\Gamma,\sigma)$. We often represent it 
as the $\bbC$-subalgebra of ${\mathcal  
B}(\ldg)$ generated by $\{T_{\gamma}|{\gamma\in\Gamma}\}$, where for  
$\gamma\in\Gamma$
\begin{eqnarray*}T_{\gamma}:\ldg & \to & \ldg,\\
\delta_\mu & \mapsto &  
\sigma(\gamma,\mu)\delta_{\gamma\mu},\end{eqnarray*}
so that an element in $\bbC(\Gamma,\sigma)$ is a finite $\bbC$-linear  
combination of the operators $T_{\gamma}$, and the convolution reads  
(for $\gamma,\mu\in\Gamma$)
$$T_{\gamma}*_{\sigma}T_{\mu}=\sigma(\gamma,\mu)T_{\gamma\mu}.$$
We shall consider several completions of $\bbC(\Gamma,\sigma)$ that we  
now explain.
The \emph{$\ell^1$-completion} (given by the norm
$\|\sum_{\gamma\in\Gamma}a_{\gamma}T_{\gamma}\|_1=\sum_{\gamma\in\Gamma} 
|a_{\gamma}|$) yields the \emph{$\ell^1$-twisted Banach algebra}  
denoted by $\ell^1(\Gamma,\sigma)$, which is the completion of  
$\bbC(\Gamma,\sigma)$ with respect to this $\ell^1$-norm. It is a  
straightforward computation to show that it is indeed a Banach algebra,  
contained in ${\mathcal B}(\ldg)$. Next we shall consider the  
\emph{operator norm}, given by
$$\|f\|_{op}=\sup\{\|f(\xi)\|_2 :  \|\xi\|_2=1\},$$
and the completion of $\bbC(\Gamma,\sigma)$ with respect to this norm  
is the \emph{twisted reduced $C*$-algebra} $C^*_r(\Gamma,\sigma)$.  
Recall that a \emph{length function on $\Gamma$} is a map  
$\ell:\Gamma\to{\bbR}_+$ satisfying:
\begin{itemize}
\item[(a)] $\ell(1)=0$, where 1 denotes the neutral element in  
$\Gamma$,
\item[(b)] $\ell(\gamma)=\ell(\gamma^{-1})$ for any $\gamma\in\Gamma$,
\item[(c)] $\ell(\gamma\mu)\leq \ell(\gamma)+\ell(\mu)$ for any
$\gamma,\mu\in\Gamma$.\end{itemize}
For $\ell$ a length function on $\Gamma$ and $s$ a positive real  
number, the \emph{$s$-weighted $\ell^2$-norm} is defined by
$$\|\sum_{\gamma\in\Gamma}a_{\gamma}T_{\gamma}\|_{s}=\sqrt{\sum_{\gamma\in\Gamma}|a_{\gamma}|^2(1+\ell(\gamma))^{2s}}$$
and the \emph{$s$-Sobolev space} is the completion of  
$\bbC(\Gamma,\sigma)$ with respect to this norm, denoted by  
$H^s_\ell (\Gamma,\sigma)$. If the length function is chosen to be 
the word length with respect to a finite set of generators for $\Gamma$,
then we just write $H^s(\Gamma,\sigma)$, omitting $\ell$ in the notation.
Finally, the \emph{space of rapidly decreasing  
functions  (with respect to the length $\ell$)} is given by
$$H^{\infty}_\ell (\Gamma,\sigma)=\bigcap_{s\geq 0}H^s_\ell (\Gamma,\sigma).$$
$H^{\infty}_\ell (\Gamma,\sigma)$ is not an algebra in general, but it is one  if  
$\Gamma$ has the Rapid Decay property (with respect to the length $\ell$), see Definition \ref{sRD}. 
In fact,  if  $\Gamma$ has the Rapid Decay property (with respect to the length $\ell$), then
$H^s_\ell (\Gamma,\sigma)$ is an algebra for $s$ large enough,
cf. Corollary \ref{remarque}. 

\begin{lemma}[Sup norm characterization]\label{supnormRD}
If $\Gamma$ has polynomial volume growth  (with respect to the length $\ell$), then the space of rapidly decreasing  
functions  (with respect to the length $\ell$) has the following 
sup norm characterization:
$$
H^{\infty}_\ell (\Gamma,\sigma)=\displaystyle\left\{ f: \Gamma \to \bbC: \displaystyle\sup_{\gamma\in\Gamma} \left(|f({\gamma})|(1+\ell(\gamma))^{s}\right)<\infty \quad \forall s\in \bbN
\right\}
$$
\end{lemma}

\begin{proof}
If $f \in H^{\infty}_\ell(\Gamma,\sigma)$, then for all $s\in \bbN$, one sees that the  function
$\gamma \mapsto |f({\gamma})|^2(1+\ell(\gamma))^{2s}$ is bounded on $\Gamma$,
therefore $\gamma \mapsto |f({\gamma})|(1+\ell(\gamma))^{s}$ is also a bounded function on $\Gamma$.

Conversely,  suppose that $ f: \Gamma \to \bbC$ is such that 
$$\sup_{\gamma\in\Gamma} \left(|f({\gamma})|(1+\ell(\gamma))^{r}\right) = C_r<\infty \quad \forall r\in \bbN.$$
Then we estimate,
$$
\sum_{\gamma\in \Gamma}|f({\gamma})|^2(1+\ell(\gamma))^{2s}\le C_r^2
\sum_{\gamma\in \Gamma} (1+\ell(\gamma))^{2(s-r)}.
$$
Since $\Gamma$ has polynomial growth (with respect to the length $\ell$), 
we see that by choosing $r$ sufficiently large,
the right hand side is finite, proving the lemma.
\end{proof}

An important step is to compute the  
$K$-theory of $C^*_r(\Gamma,\sigma)$. The following lemma shows  
that we only need to know the cohomology class of the multiplier.
\begin{lemma}\label{iso}Let $\sigma,\sigma'\in H^2(\Gamma,\Uu)$ be
two cohomologous $2$-cocycles. Then there exists an isomorphism
$$\varphi:B(\Gamma,\sigma)\to B(\Gamma,\sigma'),$$
inducing the identity map on $K$-theory. Here 
$B(\Ga,\sigma)$ is any $*$-Banach completion of $\bbC(\Ga,\sigma)$.
\end{lemma}
\begin{proof}That the cocycles $\sigma$ and $\sigma'$ are
cohomologous means that there exists $f:\Gamma\to\Uu$ such that
$\sigma'=\sigma df$, where for $\gamma_1,\gamma_2\in\Gamma$,
$df(\gamma_1,\gamma_2)=f(\gamma_1\gamma_2)f(\gamma_1)^{- 
1}f(\gamma_2)^{-1}$.
We shall define the map $\varphi:B(\Gamma,\sigma)\to
B(\Gamma,\sigma')$ on the generators
$\{T_{\gamma}\}_{\gamma\in\Gamma}$ by
$\varphi(T_{\gamma})=f(\gamma)T'_{\gamma}$, extend it by linearity to
$\bbC(\Gamma,\sigma)$ and by continuity to $B(\Gamma,\sigma)$.
Indeed, it is a *-homomorphism:
\begin{eqnarray*}\varphi(T_{\gamma}T_{\mu})&=&\sigma(\gamma,\mu)\varphi( 
T_{\gamma\mu})=\sigma(\gamma,\mu)f(\gamma\mu)T'_{\gamma\mu}\\
&=&\sigma'(\gamma,\mu)f(\gamma)f(\mu)T'_{\gamma\mu}=
f(\gamma)f(\mu)T'_{\gamma}T'_{\mu}=\varphi(T_{\gamma})\varphi(T_{\mu}),\end{eqnarray*}
bijective (it is bijective on the generators), hence induces an
isomorphism in $K$-theory, which is the identity since
$[T'_{\gamma}]=[f(\gamma)T'_{\gamma}]$ in
$K_1(B(\Gamma,\sigma'))$ for any $\gamma\in\Gamma$, the homotopy
being given by a path in $\Uu$ between $f(\gamma)$ and $1$.
Similarly at the level of $K_0$.
\end{proof}
%%%%%%%%%%%%%%%%%%%%%%%%%%%%%%%%%%%%%%%%%%%%%%%%%%%

%%%%%%%%%%%%%%%%%%%%%%%%%%%%%%

\section{Good unconditional completions}

%%%%%%%%%%%%%%%%%%%
\begin{definition}\label{good}
%{\bf Mathai, is that really needed?}
Following Lafforgue \cite{Laff}, we say that a norm $\|\ \|_{\cA}$  
on $\bbC(\Gamma,\sigma)$ is \emph{unconditional} if for any two  
elements $A=\sum_{\gamma\in\Gamma}a_\gamma T_\gamma$ and  
$B=\sum_{\gamma\in\Gamma}b_\gamma T_\gamma$ in $\bbC(\Gamma,\sigma)$,  
$|a_\gamma|\leq |b_\gamma|$ implies  
$\|A\|_{\cA}\leq\|B\|_{\cA}$. Given an unconditional norm $\|\  
\|_{\cA}$ on  $\bbC(\Gamma,\sigma)$, we denote by $\cA(\Gamma,\sigma)$  
the completion. 

For technical reasons, since we use the heat kernel approach in this paper, 
we introduce the following special case. 
Assume that an unconditional completion  $\cA(\Gamma,  
\sigma)$
of  $\C(\Gamma, \sigma)$  is such that
$$\|T_g\|_{\cA} \le C_1 e^{C_2 \ell(g)^p}, \qquad \forall g\in\Gamma, $$
for some positive constants $C_1, C_2$ independent of $g\in \Gamma$
and for some $p$ such that $1\le p<2$ which is also independent of  $g\in  
\Gamma$.
We shall call such an unconditional completion a {\em good
unconditional completion} of $\C(\Gamma, \sigma)$.
\end{definition}
\begin{remark} Note that $\ell^1$ is trivially a good unconditional  
completion, and that it is straightforward to see that the Sobolev  
completions are good unconditional completions as well. The operator  
norm is not unconditional, which means that the reduced (twisted) group $C^*$-algebra
is {\em not} an unconditional completion.
\end{remark}

\nc{\cH}{{\mathcal H}}

Let $\ell^2(\Gamma, \cH)$ denote the space of $\cH$-valued
square summable functions on the group $\Gamma$, where
$\cH$ is a separable Hilbert space with the trivial action of $\Gamma$.
  Let $\U_\cH(\Gamma, \sigma)$
denote the von Neumann algebra of all bounded linear operators
on $\ell^2(\Gamma, \cH)$ that commute with the $(\Gamma,  
\sigma)$-action.
It is a standard observation that any element  
$A\in\U_\cH(\Gamma,\sigma)$
can be represented by a strongly convergent series,
\[
A=\sum_{\gamma\in\Gamma}T_\gamma\otimes A(\gamma),
\]
where $A(\gamma)\in \B(\cH)$ is a bounded linear operator on $\cH$,  
defined by the formula
\[
A(\delta_e\otimes v)=\sum_{\gamma\in\Gamma}\delta_\gamma\otimes
A(\gamma)v, \quad v\in \cH.
\]
%
%The following lemma is elementary, cf \cite{BrSu}.
%
%\begin{lemma}\label{C}
%Let $A\in \U_\cH(\Gamma, \sigma)$,
%$A=\sum_{\gamma\in\Gamma}T_\gamma\otimes A(\gamma)$, where
%$A(\gamma)\in \B(\cH)$. Then
%\[
%\|A\| \le \sum_{\gamma\in\Gamma} \|A(\gamma)\|,
%\]
%and the right-hand side of the inequality is not necessarily  
%finite.\end{lemma}
%
One has the following useful sufficient condition, where $\otimes$ denotes 
the projective tensor product in the entire paper.
\begin{lemma}{\label D}
Let $\cA(\Gamma, \sigma)$ be a good unconditional completion of
$\C(\Gamma, \sigma)$ and $\K$ denote the algebra of compact operators
on the Hilbert space $\cH$.
  If $A\in \U_\cH(\Gamma, \sigma)$,
$A=\sum_{\gamma\in\Gamma}T_\gamma\otimes A(\gamma)$ is such that
$A(\gamma) \in \K$ and also satisfies
$
  \|A(\gamma)\|< C_5 e^{-C_6 \ell(\gamma)^2}$ for some positive constants
$C_5, C_6$, then $A \in
\cA(\Gamma, \sigma) \otimes\K$.
\end{lemma}
\begin{proof}
Observe that one has the estimate
\begin{equation}\label{e:d1}
  \# \left\{\gamma\in\Gamma\; |\;
\ell(\gamma) \le n \right\} \le C_7 e^{C_8 n},
\end{equation}
for some positive constants $C_7, C_8$, since the growth rate of
the volume of balls in $\Gamma$ is at most exponential.
We compute,
\begin{eqnarray*}
|| A ||_{\cA\otimes \K}  & = & ||\sum_{\gamma\in\Gamma}T_\gamma\otimes  
A(\gamma)||_{\cA\otimes \K} \le\sum_{\gamma\in\Gamma} ||T_\gamma||_\cA   
||A(\gamma)|| \\[+7pt]
& \le  & \sum_{\gamma\in\Gamma} C_1 e^{C_2 \ell(\gamma)^p}C_5 e^{-C_6
\ell(\gamma)^2} =  \sum_{n\in\N} \sum_{\ell(\gamma) \le n} C_1 e^{C_2  
\ell(\gamma)^p}C_5 e^{-C_6
\ell(\gamma)^2} \\[+7pt]
& \le  & \sum_{n\in\N}    C_7 e^{C_8 n} C_1 e^{C_2 n^p}C_5 e^{-C_6
n^2} < \infty.
\end{eqnarray*}
The last sum is convergent since $0\le p<2$ by the good unconditional hypothesis.
\end{proof}
In \cite{PR}, Packer and Raeburn,  inspired by A. Wasserman's thesis,  established a 
stabilization (or untwisting) trick. We will present a good unconditional version
of this, in the simple case of a discrete group $\Gamma$,  that we need in this paper.
Let $\sigma$ be a multiplier on $\Gamma$ and $\mathcal K$ be the algebra of compact  operators on
$\ell^2(\Gamma)$. Observe that for any $\Gamma$-$C^*$-algebra $A$, 
one has the following canonical isomorphism,
\begin{equation}\label{pre-PR}
\C(\Gamma, A, \sigma) \otimes {\mathcal K} \cong 
\C(\Gamma, A\otimes {\mathcal K}),
\end{equation}
where $\Gamma$ acts diagonally on the tensor product $A\otimes {\mathcal K}$, and is 
given by the given action of $\Gamma$ on $A$ and the adjoint action, 
$\gamma \mapsto {\rm Ad}(T_\gamma)$.  That is, the twisted convolution on the left hand side of \eqref{pre-PR} becomes an ordinary convolution on the right hand side of \eqref{pre-PR}: for all  $a, b \in A$, for all $U, V \in {\mathcal K}$  and for all $\gamma. \mu \in \Gamma$,
$$a \otimes {\rm Ad}(T_{\gamma}) U * b \otimes {\rm Ad}(T_{\mu}) V =a\alpha_{\gamma}(b) \otimes {\rm Ad}(T_{\gamma\mu}) UV,$$
where $\alpha$ denotes the action of $\Gamma$ on $A$. 

Recall that for any good unconditional completion  $\cA(\Gamma,  A, \sigma)$
of  $\C(\Gamma, A, \sigma)$ there are positive constants $C_1, C_2$ independent of 
$g\in \Gamma$ such that
$$\|T_g\|_{\cA} \le C_1 e^{C_2 \ell(g)^p}, \qquad \forall g\in\Gamma, $$
for some $p$ such that $1\le p<2$ which is also independent of  $g\in  
\Gamma$. But this implies that 
$$\|{\rm Ad}(T_g)\|_{\cA} \le C_1^2 e^{2C_2 \ell(g)^p}, \qquad \forall g\in\Gamma, $$
and conversely. Therefore for any good unconditional completion, one has the canonical
isomorphism
\begin{equation}\label{unconditional-PR}
\cA(\Gamma, A, \sigma) \otimes {\mathcal K} \cong 
\cA(\Gamma, A\otimes {\mathcal K}).
\end{equation}
This isomorphism is clearly also true for general unconditional completions.

%%%%%%%%%%%%%%%%%%%%%%%%%%%%%%%%%%%%%
\section{Heat kernels and the analytic twisted Baum-Connes map}\label{analytical}
%%%%%%%%%%%%%%%%%%%%%%%%%%%%%%%%%

%%%%%%%%%%%%%%%%%%%%%%%%%%%%%%%%%%%%%%%%%%%
\subsection{Spin$^\bbC$ manifolds and twisted spin$^\bbC$ Dirac operators}\label{Dirac}
%%%%%%%%%%%%%%%%%%%%%%%%%%%%%%%%%%%%%%%%%%%%%%%%%
\nc\Spin{{\rm Spin}}
\newcommand\SpinC{\operatorname{Spin}^{\bbC}}
\newcommand\spinC{\operatorname{Spin}^{\bbC}}
\nc\UU{{\rm U}}
\nc\cP{{\mathcal P}}
\nc\cF{{\mathcal F}}
\nc\cS{{\mathcal S}}

Let $M$ be a smooth $\Gamma$ manifold without boundary.
A choice of $\Gamma$-invariant 
Riemannian metric $g$ on $M$ defines a bundle of Clifford
algebras, with fibre at $z\in M$ the complexified Clifford algebra
\begin{equation}
\begin{gathered}
\Cl_z(M)=\left(\bigoplus_{k=0}^\infty (T^*_z M \otimes \bbC)^k\right)/\langle \alpha
\otimes\beta +\beta \otimes\alpha -2(\alpha ,\beta )_g,\ \alpha ,\beta
\in T^*_z M\rangle.
\end{gathered}
\end{equation}
If $\dim M=2n,$ this complexified algebra is isomorphic to the matrix algebra on
$\bbC^{2^n}.$ In particular the Clifford bundle is an associated bundle to the metric
coframe bundle, the principal $\SO(2n)$-bundle $\cF,$ where the action of
$\SO(2n)$ on the Euclidean Clifford algebra $\Cl(2n)$ is through the spin
group, which may be identified within the Clifford algebra as 
\begin{equation}
\Spin(2n)=\{v_1v_2\cdots v_{2k}\in\Cl(2n);v_i\in \bbR^{2n},\ |v_i|=1\}.
\end{equation}
The non-trivial double covering of $\SO(2n)$ is realized through the
mapping of $v$ to the reflection $R(v)\in\operatorname{O}(2n)$ in the plane
orthogonal to $v$  
\begin{equation}
p:\Spin(2n)\ni a=v_1\cdots v_{2k}\longmapsto R(v_1)\cdots R(v_{2k})=R\in\SO(2n).
\end{equation}
%Thus $\cP$ may be identified with the bundle associated to $\cF$ by the
%action of $\SO(2n)$ on $\Cl(2n)$ where
%the reflection  $R$ acts by conjugation by $a$  
%\begin{equation}
%\Cl(2n)\ni b\longmapsto aba^{-1}\in\Cl(2n).
%\end{equation}
%We therefore have a map of principal bundles 
%\begin{equation}
%\cF\longrightarrow \cP.
%\end{equation}

The $\SpinC(2n)$ group, defined as
\begin{equation}
\SpinC(2n)=\{cv_1v_2\cdots v_{2k}\in\Cl(2n);v_i\in \bbR^{2n},\ |v_i|=1,
c\in\bbC,\ |c|=1\},
\end{equation}
is a central extension of $\SO(2n),$ 
\begin{equation}
\UU(1)\longrightarrow \SpinC(2n)\longrightarrow \SO(2n),
\end{equation}
where the quotient map is consistent with the covering of $\SO(2n)$ by
$\Spin(2n),$ i.e.
\begin{equation}
\SpinC(2n) = \Spin(2n) \times_{\bbZ_2} \UU(1).
\end{equation}

The manifold $M$ is said to have a $\Gamma$-equivariant 
$\spinC$ structure, if there is an extension of the coframe bundle to a
principal $\SpinC(2n)$-bundle
\begin{equation}
\xymatrix{\UU(1)\ar[d] \ar[r]^{=}& \UU(1)\ar[d]\\
\SpinC(2n)\ar[r]\ar[d]&\cF_{L}\ar[d]\ar[r] & M\ar[d]^{||}\\\
\SO(2n)\ar[r]&\cF \ar[r]  & M,}
\label{mms3.67}\end{equation}
where $\cF_L,$ the $\SpinC(2n)$ bundle over $M$, may also be viewed as a circle bundle
over $\cF$, compatible with the $\Gamma$-action. 
The associated bundles of half spinors on $M$ are defined as
\begin{equation}
{\mathcal S}^\pm = \cF_L \times_{\SpinC(2n)} S^\pm,
\end{equation}
where $S^\pm$ are the fundamental half spin representations of $\SO(2n)$.
The $\Gamma$-invariant Levi-Civita connection determines a connection 1-form
on $\cF$, and together with the choice of a $\Gamma$-invariant connection 1-form on the circle bundle
$\cF_L$ over $\cF$, they determine a connection 1-form on the principal $\SpinC$
bundle $\cF_L$ over $M$, which is $\Gamma$-invariant.  That is, one gets a connection
\begin{equation}
\nabla^{\cS\otimes E} : C^\infty(M, {\mathcal S}^+\otimes E)  \to C^\infty(M, T^*M \otimes {\mathcal S}^+ \otimes E),
\end{equation}
defined as $\nabla^{\cS\otimes E}  = \nabla^{\mathcal S} \otimes 1 + 1 \otimes \nabla^E$, where 
$\nabla^E$ is a $\Gamma$-invariant connection on the $\Gamma$-invariant vector bundle
$E$ over $M$. Now the contraction given by Clifford multiplication defines a map
\begin{equation}
C :  C^\infty(M, T^*M \otimes {\mathcal S}^+\otimes E) \to C^\infty(M, {\mathcal S}^-\otimes E).
\end{equation}
The $\Gamma$-equivariant  
Spin$^{\mathbb C}$  Dirac operator with coefficients in $E$ is defined as the composition 
\begin{equation}
\np^{{\mathbb C}+}_E = C \circ \nabla^{\cS\otimes E}.
\end{equation}

In this section, we will define the analytic index map for an
arbitrary torsion-free discrete group $\Gamma$ and for an arbitrary
multiplier $\sigma$ on $\Gamma$ with trivial Dixmier-Douady invariant  
$\delta(\sigma)$.
%
%
%From the subsection \ref{khomology}, elements of  
%$K_0^\Gamma(\uE\Gamma)$ are equivalence classes of even dimensional   
%$\Gamma$-equivariant $K$-cycles $(M, E, \phi)$, where where 
Now let $M$ be a  manifold without boundary with a given smooth proper cocompact $\Gamma$  
action and a $\Gamma$-equivariant Spin$^{\mathbb C}$ structure, $E\to  
M$ a $\Gamma$-equivariant complex vector bundle on $M$, and $\phi :  
M \to \uE\Gamma$ a $\Gamma$-equivariant continuous map.
We will view the  $\Gamma$-equivariant  
Spin$^{\mathbb C}$  Dirac operator with coefficients in $E$ as an operator
on the $L^2$-spaces,  $\np^{\mathbb C}_E: L^2(M, {\mathcal S}^+\otimes E) \to  
L^2(M,{\mathcal S}^-\otimes E)$.

Let $c$ be an $\R$-valued $\Gamma$-equivariant  
Cech 2-cocycle on $\uE\Gamma$ and $\omega$ be a $\Gamma$-equivariant
closed 2-form on $M$ such that the $\Gamma$-equivariant
cohomology class of $\omega$ is equal to $\phi^*(c)$.
Note that $\omega$ is \emph{exact}, $\omega =  d\eta$, since $\uE\Gamma$ is  
contractible.
Define $\nabla=d+\,i\eta$.  Then $\nabla$ is a Hermitian connection
on the trivial line bundle $\cL$ over ${M}$, and the curvature
of $\nabla$ is $\ (\nabla)^2=i\,{\omega}$. Then $\nabla$
defines a projective action of $\Gamma$ on $L^2$ spinors as follows:

For $u \in L^2({M}, {\mathcal{S}}\otimes {E}\otimes \cL )$,
let $S_\gamma u  =  e^{i\varphi_\gamma}u$ (where $\varphi$ is the phase for $\sigma$ as explained in Lemma \ref{twistedIdentities}), $\;U_\gamma u  
={\gamma^{-1}}^*u$,
and $T_\gamma=U_\gamma S_\gamma$ be the composition, for all
$\gamma\in \Gamma$.  Then $T$ defines a projective  
$(\Gamma,\sigma)$-action on $L^2({M}, {\mathcal{S}}\otimes {E} \otimes \cL)$, meaning that for any  
$\gamma,\gamma'\in\Gamma$ one has
\[
T_\gamma
T_{\gamma'}=\sigma(\gamma,\gamma')T_{\gamma\gamma'}.
\]
Let $\np^{{\mathbb C}+}_{E\otimes\cL}: L^2(M, {\mathcal S}^+\otimes E\otimes \cL) \to  
L^2(M,{\mathcal S}^-\otimes E\otimes \cL)$ denote the twisted $\Gamma$-equivariant  
Spin$^{\mathbb C}$  Dirac operator.
%%%%
%%%%
\begin{lemma}
The twisted $\Gamma$-equivariant Spin$^{\mathbb C}$ Dirac operator on ${M}$,
\[
\np^{{\mathbb C}+}_{E\otimes\cL}: L^2(M, {\mathcal S}^+\otimes E\otimes \cL) \to  
L^2(M,{\mathcal S}^-\otimes E\otimes \cL),
\]
commutes with the projective $(\Gamma,\sigma)$-action.
\end{lemma}
%%%%
\begin{proof} To simplify notation, set $D_\eta = \np^{{\mathbb C}+}_{E\otimes\cL}$
and $D_0 =  \np^{{\mathbb C}+}_{E}$ where we emphasize the dependence on $\eta$. Then
$D_\eta = D_0 + ic(\eta)$, where $c(\eta)$ denotes Clifford multiplication by the one-form
$\eta$. An easy computation establishes that
$U_\gamma D_\eta = D_{{\gamma^{-1}}^*\eta} U_\gamma$ and that
$S_\gamma D_{{\gamma^{-1}}^*\eta} =   D_\eta S_\gamma \quad \text{for  
all}\;
\gamma\in\Gamma$. Then
$T_\gamma D_\eta =D_{\eta} T_\gamma$,
where $T_\gamma = U_\gamma  S_\gamma$ denotes the projective  
$(\Gamma,\sigma)$-action.
\end{proof}
%%%

\subsection{Heat kernels and the analytic index}
Recall the following well-known smoothness properties and
Gaussian off-diagonal estimates for the heat kernel, cf. \cite{Greiner,  
Ko}.
\begin{lemma}\label{E}
The Schwartz kernels $k_\pm(t,x,y)$ of the heat operators $e^{-tD^\pm  
D^\mp}$ are smooth for all $t>0$. Moreover, for any $t>0$ there are  
positive constants $C_1, C_2$ such that the following off-diagonal  
estimate holds
$$
|k_\pm (t,x,y)| \le C_1 e^{-C_2 d(x,y)^2},\quad x\in M,\quad y\in M,$$
where $d$ denotes the Riemannian distance function on $M$
\end{lemma}
For fixed $t>0$, we will use Lemma \ref{E} to show the following.
\begin{proposition}\label{off-diagonal}
Let $\cA(\Gamma, \sigma)$ be a good unconditional completion of
$\C(\Gamma, \sigma)$. Then for fixed $t>0$, the heat operators 
$e^{-tD^-D^+}$ and $e^{-tD^+D^-} $ belong to $ \cA(\Gamma, \sigma)\otimes  
{\mathcal K}_+$ and  $\cA(\Gamma, \sigma)\otimes  
{\mathcal K}_-$ respectively, where $\K_\pm$ denotes the algebra of compact operators  
on the Hilbert space $\cH_\pm= L^2(\F, \mathcal S^\pm \otimes E|_\F)$,  
and $\F$ denotes a connected fundamental
domain of the action of $\Gamma$ on $M.$
\end{proposition}
\begin{proof}
We have $e^{-tD^\pm D^\mp} \in\U_{\cH_\mp}(\Gamma, \sigma)$,  so  that
\[
e^{-tD^\pm D^\mp}
=\sum_{\gamma\in\Gamma}T_\gamma\otimes h_{t}^\pm(\gamma),
\]
where $h_{t}^\pm (\gamma)\in\B(\cH_\pm)$ has Schwartz kernel $k_\pm(t,  
x, \gamma y)$
for $x, y \in \F$.
By Lemma \ref{E}, we have $$||h_{t}^\pm (\gamma)|| \le ||k_\pm(t, x,  
\gamma y)|_\F||_\infty
\le C_1 e^{-C_2 d(\gamma)^2},$$
where $d(\gamma) = \inf\{ d(x, \gamma y): x, y \in \F\}$.

It is well known that
\begin{equation}\label{e:d}
\ell(\gamma) \le C_4 (d(\gamma)+1),
\end{equation}
for some positive constant $C_4$. From (\ref{e:d}) and
Lemma~\ref{E}, we get
\begin{equation}\label{e:2}
\|  h_{t}^\pm (\gamma)  \|\le C_5 e^{-C_6
\ell(\gamma)^2},
\end{equation}
for some positive constants $C_5, C_6$. We conclude using Lemma \ref{D}.
\end{proof}
\begin{definition}
For fixed $t>0$, define the idempotent 
$$e_t(D^+)\in  
M_2(\cA(\Gamma,\sigma)\otimes\tilde{\mathcal K})$$ as follows:
$$
e_t(D^+) = \begin{pmatrix}
e^{-tD^-D^+} & \displaystyle
e^{-{\frac{t}{2}}D^-D^+}\frac{(1-e^{-tD^-D^+})}
  {D^-D^+} D^+ \\[+11pt]
e^{-{\frac{t}{2}}D^+D^-}{D^+} &
        1- e^{-tD^+D^-}
\end{pmatrix},
$$
where $\cA(\Gamma,\sigma)\otimes
\tilde{\mathcal K}$ denotes the unital algebra associated with  
$\cA(\Gamma,\sigma)\otimes
{\mathcal K}$.
It is the analogue of the Wasserman idempotent, see e.g. Connes and  
Moscovici \cite{ConnesMosc}. Since $\cA(\Gamma,\sigma)$ is a Banach  
algebra, one has the
invariance property of $K$-theory under stable isomorphism,  
$K_\bullet(\cA(\Gamma,\sigma))
\cong K_\bullet(\cA(\Gamma,\sigma)\otimes
{\mathcal K}).
$
Using this isomorphism,  the $\cA$-{\em twisted analytic index}
is defined as
\begin{equation}\label{bc}
a\!-\!\index_\sigma^\cA(D^+)   = 
[e_t(D)] - [E_0]  \in K_0(\cA(\Gamma,\sigma)),
\end{equation}
where $t>0$ and $E_0$ is the idempotent
$$
E_0 = \begin{pmatrix} 0 & 0 \\ 0 &
        1
\end{pmatrix} \in M_2(\cA(\Gamma,\sigma)\otimes
\tilde{\mathcal K}).
$$
\end{definition}
%%%%%%%%%%
Since the difference $e_t(D^+)- E_0$ is in $M_2(\cA(\Gamma,\sigma)\otimes
{\mathcal K}), $ we see that
the right hand side of equation (\ref{bc})  is in  $K_0
(\cA(\Gamma,\sigma))$ as asserted.
%%

%%%%%%%%%%%%%%%%%%%%%%
\subsection{Topological $K$-homology and the analytic twisted Baum-Connes map}\label{khomology}
%%%%%%%%%%%%%%%%%%%%%%%%%%%%%%%%%%%%%%%%%%%%%%%%%%%%%%%
%{\bf Here say some words about spin c}. 
We shall now give a brief description of the Baum-Connes-Douglas
\cite{BC}, \cite{bd} version of the
$K$-homology groups $K_j^\Gamma(\uE\Gamma)$ ($j=0,1$).

\smallskip

The basic objects are $\Gamma$-equivariant $K$-cycles. A  
\emph{$\Gamma$-equivariant $K$-cycle} on $\uE\Gamma$ is a triple $(M,  
E, \phi)$, where:
\begin{itemize}
\item[(i)] $M$ is a  manifold without boundary with a smooth  
proper cocompact $\Gamma$-action and a $\Gamma$-equivariant
Spin$^{\mathbb C}$ structure.
\item[(ii)]$E\to M$ is a $\Gamma$-equivariant
complex vector bundle on $M$.
\item[(iii)]$\phi : M \to \uE\Gamma$ is a $\Gamma$-equivariant
continuous map.\end{itemize}

Two $\Gamma$-equivariant
$K$-cycles $(M, E, \phi)$ and
$(M', E', \phi')$ are said to be {\em isomorphic} if  there is a
$\Gamma$-equivariant
diffeomorphism $h: M \to M'$ preserving the $\Gamma$-equivariant
Spin$^{\mathbb C}$ structures on $M, M'$
such that $h^*(E') \cong E$ and
$h^*\phi' = \phi$. Let $\Pi^\Gamma(\uE\Gamma)$ denote the collection of  
all
$\Gamma$-equivariant  $K$-cycles on $\uE\Gamma$. The following  
operations on $\Gamma$-equivariant  $K$-cycles will enable us  to  
define an equivalence relation on $\Pi^\Gamma(\uE\Gamma)$.

\smallskip

\noindent {\em Bordism}: Two  $\Gamma$-equivariant $K$-cycles $(M_i,  
E_i,\phi_i) \in  \Pi^\Gamma(\uE\Gamma)$ ($i=0,1,$) are said to be  {\em  
bordant} if there
is a triple $(W, E, \phi)$, where $W$ is a
manifold with boundary $\partial W$, with a smooth proper
cocompact $\Gamma$-action and a $\Gamma$-equivariant
Spin$^{\mathbb C}$ structure;
$E\to W$ is a $\Gamma$-equivariant
complex vector bundle on $W$ and $\phi:W\to X$ is a
$\Gamma$-equivariant
continuous map such
that $(\partial W, E\big|_{\partial W}, \phi\big|_{\partial W})$ is
isomorphic to the disjoint union $(M_0, E_0, \phi_0) \cup (-M_1, E_1,
\phi_1)$. Here $-M_1$ denotes $M_1$ with the reversed
$\Gamma$-equivariant
Spin$^{\mathbb C}$ structure.

\smallskip

\noindent {\em Direct sum}: Suppose that $(M,E, \phi) \in  
\Pi^\Gamma(\uE\Gamma)$ and that
$E=E_0\oplus E_1$. Then $(M, E,\phi)$ is isomorphic to $(M, E_0,  
\phi)\cup (M, E_1,\phi)$.

\smallskip

\noindent {\em Vector bundle modification}:
Let $(M, E, \phi) \in \Pi^\Gamma(\uE\Gamma)$ and $H$ be an even  
dimensional $\Gamma$-equivariant Spin$^{\mathbb C}$
vector bundle over M. Let $\widehat M = S(H\oplus 1)$ denote the sphere
bundle of $H\oplus 1$. Then $\widehat M$ is canonically a
$\Gamma$-equivariant  Spin$^{\mathbb C}$
manifold. Let ${\mathcal S}$ denote the $\Gamma$-equivariant
bundle of spinors on $H$. Since
$H$ is even dimensional, ${\mathcal S}$ is ${\mathbb Z}_2$-graded,
$$
{\mathcal S} = {\mathcal S}^+ \oplus {\mathcal S}^-,
$$
into $\Gamma$-equivariant bundles of $1/2$-spinors on $M$.
Define $\widehat E = \pi^*({\mathcal S}^{+*} \otimes E)$, where
$\pi : \widehat M \to M$ is the projection. Finally,  set $\widehat  
\phi = \pi^*\phi$.
Then $(\widehat M,\widehat E, \widehat\phi) \in \Pi^\Gamma(\uE\Gamma)$  
is said to be obtained from
$(M, E, \phi)$ and $H$ by {\em vector bundle modification}.
%%%%%

\smallskip

Let $\;\sim\;$ denote the equivalence relation on  
$\Pi^\Gamma(\uE\Gamma)$ generated by the operations of bordism, direct
sum and vector bundle modification. Notice that $\;\sim\;$ preserves  
the parity of the dimension of
the $K$-cycle. Let
$$K_0^\Gamma(\uE\Gamma)=\Pi^\Gamma_{even}(\uE\Gamma)/\sim,$$
where
$\Pi^\Gamma_{even}(\uE\Gamma)$ denotes the set of all even dimensional
  $\Gamma$-equivariant  $K$-cycles in
$\Pi^\Gamma(\uE\Gamma)$, and let
$$K_1^\Gamma(\uE\Gamma)=\Pi^\Gamma_{odd}(\uE\Gamma)/\sim,$$
where
$\Pi^\Gamma_{odd}(\uE\Gamma)$ denotes the set of all odd dimensional
  $\Gamma$-equivariant  $K$-cycles in $\Pi^\Gamma(\uE\Gamma)$.
  
The analytic twisted Baum-Connes map is defined as
\begin{equation}\label{eq:index}
a\!-\!\mu_\sigma^{\cA} : K_0^\Gamma(\uE\Gamma) \to K_0(\cA(\Gamma,\sigma))
\end{equation}
\begin{equation}
a\!-\!\mu_\sigma^{\cA}([M,E, \phi]) =a\!-\!\index_\sigma^{\cA} (D^+),
 \end{equation}
  where $D^+$ is the twisted $\Gamma$-equivariant $\SpinC$ Dirac operator 
  defined as in section \ref{Dirac}.
  
%%%%%%%%%%%%%%%%%%%%%%%%%%%%%%%%%%%%%%%%%%%%%%%%%%%%%%%%%%%%%%%%%%%%%%%% 

%%%%%%%%%%%%%%%%%%%%%%%%%%%%%%%%%%%%%%%%%%%%%%%%%%%%%%%%%%%%%%%%%%%%%%%% 
%%%%%%%%%%%%%%%%%%%%%%%%%%%%%%%%%%%
\section{Twisted Baum-Connes conjecture in Lafforgue's settings}
%%%%%%%%%%%%%%%%%%%%%%%%%%%%%%%%%%%%%%%%%%%%%%%%%%%%%%%%%%%%%%%%%%%%%%%% 
%%%%%%%%%%%%%%%%%%%%%%%%%%%%%%%%%%%%%%%%%%%%%%%%%%%%%%%%%%
\subsection{On Lafforgue's Banach $KK$-theory}
%%%%%%%%%%%%%%%%%%%%%%%%%%%%%%%%%%%%%%
We here recall  
Lafforgue's definitions of Banach $KK$-theory in \cite{Laff}, and its  
compatibility with Kasparov $KK$-theory. A Banach algebra $A$ is called  
a $\Gamma$-Banach algebra if $\Gamma$ acts on $A$ by isometric  
automorphisms. We shall briefly sketch how Lafforgue associates to a  
pair of $\Gamma$-Banach algebras $(A,B)$ an abelian group
$$KK^{ban}_{\Gamma}(A,B).$$
An important concept in this setting is the notion of  
\emph{$\Gamma-(A,B)$-Banach bimodule} $E$: to start with, $E$ is a  
\emph{$B$-pair}, that is a pair of Banach spaces $E=(E^<,E^>)$ each of  
which is endowed with a $B$-action (left and right respectively), and with  
a $B$-valued and ${\bbC}$-linear bracket satisfying
$$\left<bx,y\right>=b\left<x,y\right>,\  
\left<x,yb\right>=\left<x,y\right>b,\  
\|\left<x,y\right>\|_B\leq\|x\|\,\|y\|,$$
(where the norms of $x$ and $y$ are taken in $E^<$ and $E^>$  
respectively). A $B$-pair $E$ is called an $(A,B)$-bimodule if it is  
endowed with a Banach algebra morphism from $A$ into ${\mathcal L}(E)$. If $A$  
and $B$ are $\Gamma$-Banach algebras, then a $B$-pair $E$ endowed with  
an isometric $\Gamma$-action is called a $\Gamma$-$B$-pair, and a  
$\Gamma-(A,B)$-Banach bimodule if $E$ is both an $(A,B)$-bimodule and a  
$\Gamma$-$B$-pair such that the morphism $A\to{\mathcal L}(E)$ is  
$\Gamma$-equivariant. Denote by $E^{ban}(A,B)$ the isomorphism classes  
of pairs $\alpha=(E,T)$, where $E$ is a ${\bbZ}_2$-graded  
$\Gamma-(A,B)$-Banach bimodule and $T\in{\mathcal L}(E)$ an operator  
reversing the graduation and such that for any $a\in A$, $[a,T]$ and  
$a(Id_E-T^2)$ are compact operator on $E$. Two cycles $\alpha$ and  
$\beta$ in $E^{ban}(A,B)$ are said \emph{homotopic} if they are the image  
of the evaluation in 0 and 1 respectively of a single element in  
$E^{ban}(A,B[0,1])$, where $B[0,1]$ denotes the Banach algebra of 
continuous maps from the interval $[0,1]$ into $B$. 
$KK^{ban}_{\Gamma}(A,B)$ is the quotient of  
$E^{ban}(A,B)$ by the equivalence relation induced by homotopy. This  
defines an abelian group, and Lafforgue's Banach $KK^{ban}$ theory is  
compatible with Kasparov's $KK$-theory in the sense that the forgetful  
morphism $\iota :E_{\Gamma}(A,B)\to E^{ban}_{\Gamma}(A,B)$ induces a  
well-defined morphism $\iota :KK_{\Gamma}(A,B)\to KK^{ban}_{\Gamma}(A,B)$  
which is functorial in $A$ and $B$ in the case where those are  
$\Gamma$-$C^*$-algebras. It is well-known that $K_*^{\Gamma}(X)\simeq  
KK_{\Gamma}(C_0(X),\bbC)$, where $X$ is any $\Gamma$-CW-complex.

\subsection{Twisted Assembly map - the idempotent method}\label{idemptent-method}
For any separable $\Gamma$-$C^*$-algebra $C$, there is  
a {\em dilation} homomorphism 
$$
\tau_{C, \Gamma} :  KK_\Gamma^{ban}(A,B) \to  KK_\Gamma^{ban}(C\otimes A,C \otimes B),
$$ 
where as in the entire paper,  $\otimes$ denotes the projective tensor
product, cf. \cite{Laff}. The following stability property is also proved in \cite{Laff}:
\begin{equation}\label{Kstability}
KK^{ban}(A,B) \cong  KK^{ban}(\K\otimes A,\K \otimes B),
\end{equation}
where $\K$ denotes the Banach algebra of compact operators on a Hilbert space, 
such as $\ell^2(\Gamma)$.
\begin{proposition}[Twisted descent map]\label{descent}
For any two $\Gamma$-$C^*$-algebras $A$ and $B$ there is a twisted descent map
\begin{equation}\label{eq:descent}
j_{\Gamma,\cA,\sigma}:KK^{ban}_{\Gamma}(A,B)\to  
KK^{ban}\left(\cA(\Gamma,A,\sigma),\cA(\Gamma,B,\sigma)\right),
\end{equation}
which is compatible with the canonical homomorphism $j_{\Gamma,\sigma}$  
of Proposition 2.1 in \cite{Ma}.
\end{proposition}
\begin{proof}
The twisted descent map is defined as the composition of the following three homomorphisms.
The first is the dilation homomorphism,
\begin{equation}\label{eq:descent1}
\tau_{\K, \Gamma} :  KK_\Gamma^{ban}(A,B) \to  KK_\Gamma^{ban}(\K\otimes A,\K \otimes B),
\end{equation}
where the action of $\Gamma$ on $\K$ is determined by $\sigma$ and is given as in  
the unconditional version of the Packer-Raeburn stabilization theorem, \eqref{unconditional-PR}. 
The second is
Lafforgue's descent homomorphism \cite{Laff},
\begin{equation}\label{eq:descent2}
j_{\Gamma, \cA} : KK_\Gamma^{ban}(\K \otimes A,\K\otimes B) \to 
KK^{ban}(\cA(\Gamma, \K \otimes A), \cA(\Gamma, \K\otimes B)),
\end{equation} 
where $\Gamma$ acts diagonally on $\K \otimes A$ and on $\K \otimes B$.
The third isomomorphism is obtained as a result of 
 the  unconditional version of the Packer-Raeburn stabilization theorem  \eqref{unconditional-PR},
together with stability of $KK^{ban}$ as above \eqref{Kstability}, 
\begin{equation}\label{eq:descent3}
KK^{ban}(\cA(\Gamma, \K \otimes A), \cA(\Gamma, \K\otimes B)) \cong
KK^{ban}\left(\cA(\Gamma,A,\sigma),\cA(\Gamma,B,\sigma)\right).
\end{equation}
The composition of the homomorphisms \eqref{eq:descent1},  \eqref{eq:descent2}, and  \eqref{eq:descent3}
yields the twisted descent map in equation \eqref{eq:descent}.
%Finally, $j_{\Gamma, \sigma}$ is well behaved under the Kasparov intersection 
%product since $j_{\Gamma}$ is.
\end{proof}

To follow Lafforgue's construction of the assembly map we shall now  
define a canonical element in $KK(\bbC,\cA(\Gamma,C_0(X),\sigma))\simeq  
K_0(\cA(\Gamma,C_0(X),\sigma))$, where $X\subset E\Gamma$ is a free cocompact  
$\Gamma$-CW-complex.
\begin{lemma}[Canonical idempotent]\label{thm:idempotent}
Take $h\in C_0(X)$ such that  
$\sum_{\gamma\in\Gamma}h(\gamma x)^2=1$ and let $\varphi$ be the phase associated to the cocycle $\sigma$. The element
\begin{equation}\label{eq:idempotent}
e(\gamma,x)=h(x)h(\gamma^{-1}x)e^{-i\varphi_{\gamma}(\gamma^{-1}x)}\in\cA(\Gamma,C_0(X),\sigma),
\end{equation}
is an idempotent, which defines a class $[e] \in K_0(\cA(\Gamma,C_0(X),\sigma))$ that is 
independent of the choice of $h$.
\end{lemma}
\begin{proof}That $e$ belongs to $\cA(\Gamma,C_0(X),\sigma)$ is clear since it is finitely supported. We now compute
\begin{eqnarray*}(e*e)(\gamma,x)&=&\sum_{g\in\Gamma}h(x)
h(g^{-1}x)e^{-i\varphi_g(g^{-1}x)}h(g^{-1}x)
h(\gamma^{-1}x)e^{-i\varphi_{g^{-1}\gamma}(\gamma^{-1}x)}\sigma(g,g^{-1}\gamma)\\
&=&h(x)h(\gamma^{-1}x)\sum_{g\in\Gamma}h(g^{-1}x)^2e^{-i(\varphi_g(g^{-1}x)+\varphi_{g^{-1}\gamma}(\gamma^{-1}x))}\sigma(g,g^{-1}\gamma)\\
&=&h(x)h(\gamma^{-1}x)\sum_{g\in\Gamma}h(g^{-1}x)^2e^{i\varphi_{\gamma}(\gamma^{-1}x)}\\
&=&h(x)h(\gamma^{-1}x) e^{i\varphi_{\gamma}(\gamma^{-1}x)} = e(\gamma, x),
\end{eqnarray*}
where the last equality follows from the relations described under Lemma \ref{twistedIdentities}.
Since the set of all $h$ as in the lemma is convex, one sees that the class $[e] \in 
K_0(\cA(\Gamma,C_0(X),\sigma))$ is independent of the choice of $h$.
\end{proof}
%
%Now $e$ defines a class $[e] \in K_0(\cA(\Gamma,C_0(X),\sigma))$.
We denote by  
$$p:KK^{ban}\left(\cA(\Gamma,C_0(X),\sigma),\cA(\Gamma,\sigma)\right)\to  
K_0(\cA(\Gamma,\sigma)),$$ the map determined by the idempotent $e$, 
i.e. $p(\xi) = [e] \otimes_{\cA(\Gamma,C_0(X),\sigma)} \xi \in K_0(\cA(\Gamma,\sigma))$
for all $\xi \in  KK^{ban}\left(\cA(\Gamma,C_0(X),\sigma),\cA(\Gamma,\sigma\right)$, cf. 
Lemma \ref{eq:idempotent}, as done by  
Lafforgue in \cite{Laff} page 42 in the untwisted case.

 The {\em twisted assembly map},
\begin{equation}\label{eq:iassembly}
t\!-\!\mu_{\sigma}^{\cA}:K_*^{\Gamma}(\uE\Gamma)\simeq KK_{\Gamma}(C_0(\uE\Gamma),\bbC)\to  
K_*(\cA(\Gamma,\sigma)),
\end{equation}
 is then  
defined as the inductive limit over cocompact $\Gamma$-CW-complexes $X$  
of the following maps:
$$t\!-\!\mu_{\sigma}^{\cA, X}:K_*^{\Gamma}(X)\simeq KK_{\Gamma}(C_0(X),\bbC)\to  
K_*(\cA(\Gamma,\sigma)),$$
where each map $t\!-\!\mu_{\sigma}^{\cA, X}$ is given as the composition
$p\circ j_{\Gamma,\cA,\sigma}\circ\iota$,  
that is,  
{\small
$$KK_{\Gamma}(C_0(X),\bbC)\stackrel{\iota}{\to} KK^{ban}_{\Gamma}(C_0(X),\bbC) 
)\stackrel{ j_{\Gamma,\cA,\sigma}}{\to} KK^{ban}\left(\cA(\Gamma,C_0(X),\sigma),\cA(\Gamma,\sigma)\right)
\stackrel{p}{\to}  
K_*(\cA(\Gamma,\sigma)).$$}

%%%%%%%%%%%%%%%%%%%%%%%%%%%%%%%%%%%%%%
\subsection{On the equivalence of the analytic twisted Baum-Connes and  twisted assembly maps}
\label{equivalence}
%%%%%%%%%%%%%%%%%%%%%%%%%%%%%%%%%%%%%%%

We sketch the equivalence of the twisted assembly maps given by 
equations \eqref{eq:iassembly} and \eqref{eq:index}.	

As in section \ref{khomology}, let $(M, E, \phi)$ denote a $\Gamma$ equivariant $K$-cycle. 
Then the analytic twisted Baum-Connes map is defined as
 in \eqref{eq:index} in terms of the analytic index, 
\begin{equation}
a\!-\!\mu_\sigma^{\cA} : K_0^\Gamma(\uE\Gamma) \to K_0(\cA(\Gamma,\sigma))
\end{equation}
\begin{equation}
a\!-\!\mu_\sigma^{\cA}([M,E, \phi]) = a\!-\!\index_\sigma^{\cA} (D^+),
 \end{equation}
  where $D$ is the twisted $\SpinC$ Dirac operator defined as in section \ref{Dirac}.
  
  On the other hand, section \ref{idemptent-method} defines a twisted assembly
  map in terms of the class of an idempotent, $[e] \in K_0(\cA(\Gamma,C_0(X),\sigma))$,
  as 
\begin{equation}
t\!-\!\mu_\sigma^{\cA} : K_0^\Gamma(\uE\Gamma) \to K_0(\cA(\Gamma,\sigma))
\end{equation}
\begin{equation}
t\!-\!\mu_\sigma^{\cA}([M,E, \phi]) = t\!-\!\index_\sigma^{\cA} (D^+),
 \end{equation}
where $ t\!-\!\index_\sigma^{\cA} (D^+) =  [e] \otimes_{\cA(\Gamma,C_0(X),\sigma)} 
 j_{\Gamma,\cA,\sigma}([M, E, \phi]) \in K_0(\cA(\Gamma,\sigma))$. 
 
A direct application of the scheme of section 4, \cite{Ka2}, 
establishes the following index theorem,
\begin{equation}\label{a-index=t-index}
 a\!-\!\index_\sigma^{\cA} (D^+)  =  t\!-\!\index_\sigma^{\cA} (D^+)  \in K_0(\cA(\Gamma,\sigma)).
\end{equation}

 Therefore the 	analytic twisted Baum-Connes map $a\!-\!\mu_\sigma^{\cA}$ 
 and 	the  twisted assembly map $t\!-\!\mu_\sigma^{\cA}$ are equal, so we will
henceforth denote either of these by $\mu_\sigma^{\cA}$.

%%%%%%%%%%%%%%%%%%%%%%%%%%%%%%%%%%%%%%
\subsection{Unconditional analog of the twisted Baum-Connes conjecture} 
%%%%%%%%%%%%%%%%%%%%%%%%%%%%%%%%%%%%%%%
The following conjecture is natural in view of the above computations  
combined with Lafforgue's work, and it amounts to a twisted Bost  
conjecture in case where we choose the unconditional completion to be  
$\ell^1$.
\begin{conjecture}\label{TwistedBost} Let $\Gamma$ be a countable group  
and $\sigma$ a
multiplier on $\Gamma$ with
trivial Dixmier-Douady invariant. Then for any unconditional completion
$\cA(\Gamma,\sigma)$  of
$\C(\Gamma,\sigma),$ the twisted assembly map
%%%
$$
\mu_\sigma^{\cA} : K_j^\Gamma( \uE\Gamma) \rightarrow K_j
(\cA(\Gamma,
\sigma)),\qquad j=0,1,
$$
%%%
is an isomorphism.
\end{conjecture}
%%%
This conjecture is strongly related to a twisted version of the  
Baum-Connes conjecture (see \cite{BC}).
\begin{conjecture}\label{TwistedBC} Let $\Gamma$ be a countable group  
and $\sigma$ a
multiplier on $\Gamma$ with
trivial Dixmier-Douady invariant. Then the twisted assembly map
%%%
$$
\mu_\sigma : K_j^\Gamma( \uE\Gamma) \rightarrow  
K_j(C^*_r(\Gamma,\sigma)),\qquad j=0,1,
$$
%%%
is an isomorphism.\end{conjecture}
%%%
To prove Conjecture \ref{TwistedBC} in some cases, we first prove  
Conjecture \ref{TwistedBost} and deduce Conjecture \ref{TwistedBC} from  
Conjecture \ref{TwistedBost} when the groups in addition have property  
RD, using Proposition \ref{KtheoryIso}. To prove Conjecture~\ref{TwistedBost} in case where the group $\Gamma$ is in Lafforgue's  
class ${\mathcal C}'$ we need to first recall some facts and definitions.
%%%
Let $A$  be a proper $\Gamma$-C*-algebra. Then a Dirac element $\alpha\in  
KK^\Gamma_0(A,{\bbC})$ and a dual Dirac element $\beta\in KK^\Gamma_0({\bbC},A)$ 
satisfy the following conditions,
\begin{equation}
\begin{array}{lcl}
\alpha \otimes_{\bbC} \beta & = & 1 \in KK^\Gamma_0(A, A)\\[+7pt]
\beta \otimes_{A} \alpha & = & \gamma \in KK^\Gamma_0(\bbC, \bbC),
\end{array}
\end{equation}
where $\gamma$ is the idempotent as defined by Kasparov in \cite{Kasparov}, Lafforgue in  
\cite{Laff} or Valette in \cite{Valette}. The Dirac element $\alpha$ gets its name as it 
is constructed using a $\SpinC$ Dirac operator.

\begin{definition} We say that a group $\Gamma$ has the \emph{Banach  
Dirac-Dual Dirac property} if the element $\gamma\in KK_\Gamma  
(\bbC,\bbC)$ is trivial in $KK^{
ban}_\Gamma (\bbC,\bbC)$.\end{definition}
%%%%%%
Recall that Lafforgue's class ${\mathcal C}'$ defined in \cite{Laff}  
contains all countable discrete groups acting properly and by  
isometries either on a Hilbert space (those are said to have the  
Haagerup property, see \cite{les_welches}, which include amenable groups,  
free groups and the property is closed under free and direct products),  
on a strongly bolic space (e.g. CAT(0) groups, hyperbolic groups due to  
Mineyev and Yu \cite{Mineyev_Yu}), or on some non-positively curved  
Riemannian manifolds (as linear groups).
%%%
\begin{theorem}[Lafforgue \cite{Laff}] Any group $\Gamma$ in the class  
${\mathcal C}'$ has the Banach Dirac-Dual Dirac property.\end{theorem}
%%%%%
\begin{theorem}\label{main1} Suppose that $\Gamma$ is a
discrete group that has the {Banach Dirac-Dual Dirac property}, and
that $\sigma$ is a multiplier on $\Gamma$ with
trivial Dixmier-Douady invariant. Then Conjecture~\ref{TwistedBost}  is  
true.

If in addition, $\Gamma$ has property RD, then  
Conjecture~\ref{TwistedBC}  is true.
\end{theorem}
\begin{proof}[Sketch] It follows from Lafforgue's work that we can find  
a proper $\Gamma$-C*-algebra $A$ and a Dirac element $\alpha\in  
KK^\Gamma_i(A,{\bf C})$ and a dual Dirac element $\beta\in KK^\Gamma_i({\bf C},A)$ 
such that
$$\gamma=\beta\otimes_A\alpha=1\ \ \ \hbox{in}\ KK_\Gamma^{ban}({\bf  
C},{\bf C}).$$
Then consider the following commutative diagram:
$$\xymatrix{
K^\Gamma_*(\underline{E}\Gamma)\ar[rr]^-{\otimes_{\bf  
C}\beta}\ar[d]^{\mu_{\sigma}^{\cA}} & &  
KK^\Gamma_*(\underline{E}\Gamma,A)\ar[d]^{\mu_*^{\Gamma,A}}_{\simeq}\ar[rr]^-{\otimes_A\alpha} & &  
K^\Gamma_*(\underline{E}\Gamma)\ar[d]^{\mu_{\sigma}^{\cA}}\\
K_*(\cA(\Gamma,\sigma))\ar[rr]_-{\bigotimes j_\Gamma(\beta)} & &  
K_*(\cA(\Gamma,A,\sigma))\ar[rr]_-{\bigotimes j_\Gamma(\alpha)} & &  
K_*(\cA(\Gamma,\sigma)),
}$$
with the fact that, using Proposition \ref{descent}, composites on the top  
and the bottom lines are identity.
\end{proof}

%%%%%%%%%%%%%%%%%%%%%%%%%%%%%%%%%%%%%%%%%%%%%%%%%%%%%%%%%%%%%%%%%%%%%%%% 
%%%%%%%%%%%%%%%%%%%%%%%%%
\section{On the range of the trace, conjectures and applications}
%%%%%%%%%%%%%%%%%%%%%%%%%%%%%%%%%%%%%%%%%%%%%%%%%%%%%%%%%%%%%%%%%
\subsection{Characteristic classes}
%%%%%%%%%%%%%%%%%%%%%%%%%%%%%%%%%%%%%%%%%%%%%%%%%%%%%%%%%%%%%%%%%%%%%%%% 
%%%%%%%%%%%%%%%%%%%%%%%%%
We recall some basic facts about some well-known
characteristic classes that will be used in this paper,
cf. \cite{Hir}.

Let $E\to M$ be a Hermitian vector bundle over the compact
manifold $M$ that has dimension $n=2m$.
   The {\em Chern classes} of $E$, $c_j(E)$, are by
definition {\em integral } cohomology classes.
The {\em Chern character} of $E$, $\Ch(E)$, is a rational
cohomology class
$$
\Ch(E)=\sum_{r=0}^{m} \Ch_r(E),
$$
where $\Ch_r(E)$ denotes the component of $\Ch(E)$ of degree
$2r$. Then $\Ch_0(E) = \rank(E)$, $\Ch_1(E) = c_1(E)$
and in general
$$
\Ch_r(E) = \frac{1}{r!} P_r(E)
\in H^{2r}(M, \mathbb Q),
$$
where $P_r(E) \in  H^{2r}(M, \mathbb Z)$
is a polynomial in the Chern classes of degree
less than or equal to $r$ with {\em integral} coefficients, that is
determined inductively by the Newton formula
$$
P_r(E) - c_1(E) P_{r-1}(E)+\ldots +
(-1)^{r-1}c_{r-1}(E) P_1(E)
+(-1)^{r} r c_{r}(E) =0,
$$
and by $P_0(E) = \rank(E)$. The next two terms are
$P_1(E) = c_1(E)$, $P_2(E) = c_1(E)^2 - 2 c_2(E)$.

The Todd-genus characteristic class of the
Hermitian vector bundle $E$ is a rational cohomology class in
$H^{2\bullet}(M,\mathbb Q)$,
$$
\Todd(E) = \sum_{r=0}^{m}\Todd_r(E),
$$
where $\Todd_r(E) $ denotes the component of
$\Todd (E)$ of degree $2r$.  Then $\Todd_r(E)  = B_r Q_r(E)$, where
$Q_r(E)$  is a polynomial in the Chern classes of
degree less than or equal to $r$, with
{\em integral} coefficients, and $B_r \ne 0, B_r \in \mathbb Q$ are the
Bernoulli numbers. In particular, $\Todd_0(E) = B_0 Q_0 = 1$.

{\em For the rest of this section, we use the notation of Section \ref{Dirac}}.

%%%%%%%%%%%%%%%%%%%%%%%%%%%%
\subsection{An $L^2$ index theorem}\label{l2index}
%%%%%%%%%%%%%%%%%%%%%%%%%%%%%

Let  $\;\tau\;$ be the canonical trace on $\cA(\Gamma,\sigma)\;$
defined by evaluation at the identity element of $\Gamma$. It
induces a linear map
$$
[\tau] : K_0 (\cA(\Gamma,
\sigma)) \to \mathbb R,
$$
which is called the {\em trace map} in $K$-theory.
Explicitly, first $\;\tau\;$ extends to matrices with
entries in $\cA(\Gamma, \sigma)\;$ as (with Trace denoting
matrix trace):
\[
     \tau(f\otimes r) = {\mbox{Trace}}(r)
\tau(f).
\]

Then the extension of $\;\tau\;$ to $K_0$ is given by
$\;[\tau]([e]-[f]) = \tau(e) - \tau(f),\;$ where $e$ and $f$
are idempotent matrices with entries in $\cA(\Gamma,\sigma)$.

We will compute $[\tau]\circ \mu_\sigma^\cA(K_0^\Gamma(\uE\Gamma))$
as follows.
$$
\begin{array}{rcl}
[\tau]\circ \mu_\sigma^\cA ([M,E, \phi]) & = & [\tau](   [ e_t(D) ] -  
[E_0])\\[+7pt]
& = & \tau(e^{-tD^-D^+} ) - \tau(e^{-tD^+D^-} ) \\[+7pt]
& = &  c_0 \displaystyle \int_{M/\Gamma} \Todd
(M)\wedge e^{\omega}
\wedge\Ch(E),
\end{array}
$$
where $D^+ = \np^{{\mathbb C}+}_{E\otimes\cL}$ denotes the twisted 
$\Gamma$-equivariant Dirac operator and $D^-$ denotes its adjoint.
Here the local index theorem is used to deduce the last line,
cf. the Appendix in  \cite{Ma2}.
Here $c_0= {1}/{(2\pi)^{n/2}}$ is the universal constant determined by 
the Atyiah-Singer index theorem, see \cite{AtSi} ,
$n = \dim M$, $\Todd$ and $\Ch$ denote
the Todd-genus and the Chern character respectively, $\omega$ is the 
curvature of the connection on the trivial line bundle $\cL$ 
that is described in Section \ref{Dirac}. 
This theorem is also a consequence of section \ref{equivalence}.

%%%%%%%%%%%%%%%%%%%%%%%%%%%%%%%%%%%%%%%%%%%%%%%%%%%%%%%%%%%%%%%%%%%%
\subsection{Range of the canonical trace}\label{sect:range}
%%%%%%%%%%%%%%%%%%%%%%%%%%%%%%%%%%%%%%%%%%%%%%%%%%%%%%%%%%%%%%%%%%%
Here we will present some consequences  of the twisted Bost conjecture
above and the twisted $L^2$ index theorem
described in subsection \ref{l2index} above.

The following result is an easy modification of a result in  \cite{Ma}.
%%%
\begin{theorem}[Range of the trace theorem]\label{range} Suppose that
$(\Gamma, \sigma)$ satisfies Conjecture~\ref{TwistedBost} .
Then the range of the canonical trace on $K_0(\cA(\Gamma, \sigma))$ is given by
$$
\left\{ c_0 \int_{M/\Gamma}\Todd
(M)\wedge e^{\omega}
\wedge\Ch(E) ;\; \mbox{\rm for all}\;(M, E,
\phi)
\in \Pi_{\rm even}^\Gamma(\uE\Gamma)\right\}.
$$
\end{theorem}
%%%%
\begin{remarks}
The set
$$\left\{ c_0 \int_{M/\Gamma} \Todd (M)\wedge e^{\omega}
\wedge\Ch(E) ;
\;\mbox{\rm for all}\; (M, E, \phi)
\in \Pi_{\rm even}^\Gamma(\uE\Gamma)\right\},$$ is a
countable discrete subgroup of $\R$, but  it is not in general a  
subgroup of $\Z$.
\end{remarks}
%%%
\begin{remarks}
When $\Gamma$ is the fundamental group of a compact Riemann surface of  
positive
genus, it follows from \cite{Rieff} in the genus one case,  and  
\cite{CHMM}
in the general case, that the set
$$\left\{ c_0 \int_{M/\Gamma} \Todd (M)\wedge
e^{\omega}
\wedge\Ch(E) ;  \;\mbox{\rm for all}\; (M, E, \phi)
\in
  \Pi_{\rm even}^\Gamma(\uE\Gamma)\right\},$$
reduces to the countable discrete group $\Z + \theta \Z$, where
$\theta \in [0,1)$ corresponds to the multiplier $\sigma$ under the  
isomorphism
$H^2(\Gamma; {\rm\bf U}(1))\cong \R/\Z$.
\end{remarks}
%%%
\begin{proof}[Proof of Theorem {\rm {\ref{range}}}]
By hypothesis, the twisted assemby map $\mu_\sigma^\cA$ is
an isomorphism. Therefore to compute the
range of the trace map on
$K_0(\cA(\Gamma, \sigma))$, it suffices to compute the range of
the trace map on elements of the form
$$\mu_\sigma^\cA([M, E, \phi]), \qquad
[M, E, \phi] \in K_0^\Gamma(\uE\Gamma).$$
Here $(M, E, \phi)\in \Pi_{\rm even}^\Gamma(\uE\Gamma)$.
By the $L^2$ index theorem described in section \ref{l2index} above,  
one has
\begin{equation*}
     [\tau](\mu_\sigma^\cA([M, E, \phi])) = c_0 \int_{M/\Gamma} \Todd  
(M)\wedge e^{\omega}\wedge\Ch(E),
\end{equation*}
as desired.
\end{proof}
Therefore we deduce the following.
%%%%%%%%%%%%%%
\begin{cor} Suppose that
$(\Gamma, \sigma)$ satisfies Conjecture \ref{TwistedBC}, then the range  
of the trace map on $K_0(C^*_r(\Gamma, \sigma))$ is
%%%
$$\left\{ c_0 \int_{M/\Gamma}\Todd
(M)\wedge e^{\omega}
\wedge\Ch(E) ;\; \text{for all}\;(M, E,
\phi)
\in \Pi_{even}^\Gamma(\uE\Gamma)\right\}.
$$
%%%
\end{cor}
%%%%%%%%%%%%%%%%%%%%%%%%%%%%%%%%%%%%%%%%%%%
\subsection{The 3 and 4 dimensional cases}\label{3&4}
%%%%%%%%%%%%%%%%%%%%%%%%%%%%%%%%%%%%%%%%%%%%
We explicitly determine the range of the trace in the special case when $\Ga$ is torsion-free and $B\Gamma$  is either a three or a four dimensional smooth compact manifold. In the three dimensional case we get the following.
\begin{theorem}\label{3D} Let $\Ga$ be a torsion-free group such that $(\Gamma, \sigma)$ satisfies Conjecture \ref{TwistedBost}, and such that $B\Gamma$ is a smooth, compact oriented three dimensional manifold.
Then the range of the trace map is
\begin{equation}\label{eq:3D}
[\tr] (K_0 (\cA(\Gamma, \sigma)) ) = \Z + \sum_{i=1}^{b_1} \Z \theta_i,
\end{equation}
where for $i=1,\dots, b_1$, $\theta_i= c_0 \left<\eta_i\cup\omega,[B\Ga]\right>$ and the $\eta_i$'s are generators for $H^1(B\Ga,\Z) \cap H^1(B\Gamma, \bbR) $,  $b_1 = \dim H^1(B\Gamma, \bbR) $ and $\sigma=e^{\omega}$.\end{theorem}
\begin{proof}
Since any smooth, compact oriented three dimensional manifold is a spin manifold, 
it satisfies Poincar\' e duality, 
$$
K_0(B\Gamma) \cong K^1(B\Gamma).
$$
The range of the trace $[\tr] (K_0 (\cA(\Gamma, \sigma)) ) $   simplifies to  
\begin{equation}\label{eq:3D1}
\left\{ c_0 \int_{B\Gamma} \Todd (B\Gamma)\wedge e^{\omega} \wedge\Ch^{odd}(E) ;\; \text{for all}\; E
\in K^1(B\Gamma) \right\}.
\end{equation}
For dimension reasons, $\Todd (B\Gamma) = 1$, $e^{\omega} = 1 + \omega$ and 
$\Ch^{odd}(E) = c^{odd}_1(E) + \Ch^{odd}_3(E)$.  Therefore equation \eqref{eq:3D1} reduces to 
$$
\left\{ c_0 \int_{B\Gamma} c^{odd}_1(E) \wedge \omega + c_0 \int_{B\Gamma} \Ch^{odd}_3 (E)  ;\; \text{for all}\; E
\in K^1(B\Gamma) \right\}.
$$
By the Atiyah-Singer index theorem \cite{AtSi}, one knows that 
$$
c_0 \int_{B\Gamma} \Ch^{odd}_3 (E) \in \mathbb Z \qquad \text{for all}\; E
\in K^1(B\Gamma).
$$
The proof is concluded from the fact that $c^{odd}_1(E) = c^{odd}_1(\det E) \in H^1(B\Ga,\Z) \cap H^1(B\Gamma, \bbR)$. 
\end{proof}

We now turn to the four dimensional case. Let $Q(a,b) = \langle a\cup b, [B\Gamma]\rangle$, for $a,b \in  
H^2(\Gamma, \R)$, 
be the intersection form on $B\Gamma$. Define the
linear functional $T_\omega : H^2(\Gamma, \Z) \to \R$ as $T_\omega
(a) = Q(\omega,a)$. Then the following is a consequence of
Theorem \ref{range} and the proof of Theorem 2.5 in \cite{MM}.
%%%%
\begin{theorem} Let $\Ga$ be a torsion-free group such that $(\Gamma, \sigma)$ satisfies Conjecture \ref{TwistedBost}, and such that $B\Gamma$ is a smooth, compact oriented four dimensional manifold.
Then the range of the trace map is
\begin{equation}\label{4D}
[\tr] (K_0 (\cA(\Gamma, \sigma)) ) = \Z\theta + \Z + B,
\end{equation}
where $2(2\pi)^2\theta = \langle[\omega\cup \omega], [B\Gamma]\rangle$,
and $B= {\rm range}(T_\omega)$.
\end{theorem}
\begin{remarks} Here $\omega$ is as in subsection \ref{l2index}.
If $a_1, \dots, a_{r}$ are generators of $H^2(B\Gamma, \Z) \cap H^2(B\Gamma, \bbR) $,  
where $r = \dim H^2(B\Gamma, \bbR) $, then we can
express equation (\ref{4D}) as,
$$
[\tr] (K_0 (\cA(\Gamma, \sigma)) ) = \Z\theta + \Z + \sum_{j=1}^{r} 
\Z\theta_j,
$$
where $\theta_j =  \langle  \omega\cup a_j, [B\Gamma]\rangle$ for $j=1, \dots, r$.
\end{remarks}

%%%%%%%%%%%%%%%%%%%%%%%%%%%%%%%%%%%%%%%%%%%
\subsection{The trace conjecture for unconditional twisted group completions}
%%%%%%%%%%%%%%%%%%%%%%%%%%%%%%%%%%%%%%%%%%%%
The calculations done earlier in the section validate the following bold conjecture.

\begin{conjecture}\label{conj:trace}
Let $\Ga$ be a torsion-free group such that $(\Gamma, \sigma)$ satisfies Conjecture \ref{TwistedBost}, and such that the classifying space $B\Gamma$ is a smooth, compact,  oriented manifold.\\

\noindent{\rm (1) (Even dimensional case)}  Suppose that $B\Gamma$ is of dimension $2n$.
If $a_1(j), \dots, a_{b_{2j}}(j)$ are generators of $H^{2j}(B\Gamma, \Z) \cap H^{2j}(B\Gamma, \bbR) $,  
where $b_{2j} = \dim H^{2j}(B\Gamma, \bbR) $, then  the range of the trace map is
\begin{equation}\label{evenD}
[\tr] (K_0 (\cA(\Gamma, \sigma)) ) = \Z + \Z\theta  + \sum_{j=1}^{n-1}\sum_{k=1}^{b_{2j}}  \Z r_{k, j, n} \theta_k(j),
\end{equation}
where $\theta_k(j) =  \langle [ \omega^{n-j}\cup a_k(j)], [B\Gamma]\rangle$ for $k=1, \dots, b_{2j}$,
$2(2\pi)^n\theta = \langle[\omega^n], [B\Gamma]\rangle$ and $r_{k, j, n} $ are universal constants.\\

\noindent{\rm (2) (Odd dimensional case)}
Suppose that $B\Gamma$ is  of dimension $2n-1$.
If $a_1(j), \dots, a_{b_{2j-1}}(j)$ are generators of $H^{2j-1}(B\Gamma, \Z) \cap H^{2j-1}(B\Gamma, \bbR) $,  
where $b_{2j-1} = \dim H^{2j-1}(B\Gamma, \bbR) $, then  the range of the trace map is
\begin{equation}\label{oddD}
[\tr] (K_0 (\cA(\Gamma, \sigma)) ) = \Z +  
%\frac{1}{(n-1)!} \sum_{k=1}^{b_1} \Z  \theta_k(1)+  
% \sum_{k=1}^{b_{2n-3} }\Z  \theta_k(n-1)  +  
\sum_{j=1}^{n-1} 
\sum_{k=1}^{b_{2j-1}}  \Z r_{k, j, n}' \theta_k(j),
\end{equation}
where $\theta_k(j) =  \langle [ \omega^{n-j}\cup a_k(j)], [B\Gamma]\rangle$ for $k=1, \dots, b_{2j-1}$
and $r_{k, j, n} '$ are universal constants.
\end{conjecture}

\begin{remark}
In section \ref{3&4}, we have a more explicit form of Conjecture \ref{conj:trace},
whenever the dimension of  $B\Gamma$ is less than or equal to 4.
\end{remark}

%By the earlier calculations, $r_{k, j, n}' \equiv 1$ whenever $n \le 2$.

%%%%%%%%%%%%%%%%%%%%%%%%%%%%%%%%%%%%%%%%%%%%%%%%%%%%%%%%%%%%%%%%%%%%%%%% 
%%%%%%
\section{Generators for $K_1$ of twisted group algebra completions of surface groups}
%%%%%%%%%%%%%%%%%%%%%%%%%%%%%%%%%%%%%%%%%%%%%%%%%%%%%%%%%%%%%%%%%%%%%%%% 
%%%%%%%%%%%%%%%%%%%%%%%%%%%%%%%%%%%
In this section we shall focus on the degree one part of conjecture  
\ref{TwistedBost}, and more precisely we shall identify the generators  
for $K_1$ of a fundamental group of a Riemann surface. We assume that  
all the groups are torsion-free. We shall more generally obtain partial  
results in low-dimensional homology ($B\Gamma$ of dimension  
less than or equal to 2). In this case, there exist (as shown by  
Natsume in \cite{Natsume} in the untwisted case) natural homomorphisms
\begin{eqnarray*}
\beta_t:{H_1(\Gamma,{\Z})}=\Gamma^{\rm ab} & \to & {K_1(B\Gamma)}\\
\beta_a^\sigma: {H_1(\Gamma,{\Z})}=\Gamma^{\rm ab} & \to &
{K_1(\cA(\Gamma, \sigma))}
\end{eqnarray*}
such that $\beta_a^\sigma =\mu_1^\sigma\circ\beta_t$. 
Here $\cA(\Gamma, \sigma)$ is any good unconditional completion 
of the twisted group algebra $\bbC(\Gamma, \sigma)$. Defining  
$\beta_t:\Gamma^{\rm ab}\to
K_1(B\Gamma)$ has been done by Valette in \cite{Valette}, and we
recall here the construction. Since $\pi_1(B\Gamma)=\Gamma$, an
element $\gamma\in\Gamma$ can be viewed as a pointed continuous map
$\gamma:S^1\to B\Gamma$, inducing a map in $K$-homology,
$$\gamma_*:K_1(S^1)\to K_1(B\Gamma).$$
The generator of $K_1(S^1)\simeq{\bbZ}$ can be described by the class
of the cycle $(\pi,D)$ where $\pi$ is the representation of $C(S^1)$ on
$L^2(S^1)$ by pointwise multiplication and
$$D=-i\frac{t}{dt}.$$
An element $\gamma\in\Gamma$ gets then mapped to the class of the
cycle $\gamma_*(\pi,D)=(\gamma_*\pi,D)$, where for $X$ a compact
subset of $B\Gamma$ containing $\gamma(S^1)$ and $f\in C(X)$,
$\gamma_*\pi(f)=\pi(f\circ\gamma)$ is the pointwise multiplication by
$f\circ\gamma$ on $L^2(S^1)$. In other terms one defines
\begin{eqnarray*}\tilde{\beta}_t:\Gamma & \to & K_1(B\Gamma)\\
\gamma & \mapsto & [(\gamma_*\pi,D)].\end{eqnarray*}
According to Valette in \cite{Valette}, the map
$\tilde{\beta}_t:\Gamma\to K_1(B\Gamma)$ is a group homomorphism and
hence factors through
$$\beta_t:\Gamma^{\rm ab}\to K_1(B\Gamma).$$
To define $\beta_a^\sigma:\Gamma^{\rm ab}\to{K_1(\cA(\Gamma,
\sigma))}$, simply map a representative $[\gamma]$ of $\Gamma^{\rm ab}$
to the class
$[T_{[\gamma]}]$ of the invertible operator $T_{[\gamma]}$ in
$\cA(\Gamma, \sigma)$. The following is then an adaptation of a result due  
to Natsume \cite{Natsume}, and the proof we give here is an easy twist  
of the one given in \cite{Valette}, taken from \cite{Bettaieb-Matthey}.  
It explains how the map $\beta_t$  is related to the maps
$\mu_1^{\sigma}$ and $\beta_a^{\sigma}$.
\begin{proposition}\label{K1twist} For $\sigma$ a multiplier on  
$\Gamma$ with
$\delta(\sigma)=0$, we have that
$\beta_a^{\sigma}=\mu_{\sigma}^{\cA}\circ\beta_t$.
\end{proposition}
\begin{proof}It is enough to see that
$$\tilde{\beta}_a^{\sigma}=\mu_{\sigma}^{\cA}\circ\tilde{\beta}_t:\Gamma 
\to
K_1(\cA(\Gamma,\sigma)).$$
For $\gamma\in\Gamma$, denote by $\gamma$ the (unique) homomorphism
${\bbZ}\to\Gamma$ such that $\gamma(1)=\gamma$. Consider then the
diagram
$$\xymatrix{
{\bbZ}\ar[rrr]_{\gamma}\ar[rd]^{\beta_a^{\sigma}}\ar[dd]_{\beta_t} &
& &
{\Gamma}\ar[dl]^{\tilde{\beta}_a^{\sigma}}\ar[dd]^{\tilde{\beta}_t}\\
   & {K_1(\cA({\bbZ},\sigma))}\ar[r]_{\gamma_*} &  
{K_1(\cA(\Gamma,\sigma))} & \\
{K_1(S^1)}\ar[ur]_{\mu_{\sigma}^{\cA}}
\ar[rrr]_{\gamma_*}
& & & {K_1(B\Gamma)}\ar[ul]_{\mu_{\sigma}^{\cA}}
}$$
where by abuse of notation $\sigma$ denotes also the multiplier
$\sigma$ restricted to ${\bbZ}$. That
$\tilde{\beta}_a^{\sigma}\circ\gamma=\gamma_*\circ\beta_a^{\sigma}$
is a simple computation,
$\tilde{\beta}_t\circ\gamma=\gamma_*\circ\beta_t$ by definition of
$\tilde{\beta}_t$, and
$\gamma_*\circ\mu_1^{\sigma}=\mu_1^{\sigma}\circ\gamma_*$ by
naturality of the twisted assembly map (see \cite{Ma}). That
$\beta_a=\mu_1^{\bbZ}\circ\beta_t$ follows from the proof of the
isomorphism \ref{TwistedBC} for ${\bbZ}$ (see \cite{Ma}) and
we conclude the proof by a diagram chase.\end{proof}
\begin{cor}\label{thm:rs} Let $\Gamma_g$ be the fundamental group of a 
compact Riemann surface of genus $g\ge 1$.  
Then the map $\beta_a^\sigma$ is an isomorphism.
\end{cor}
\begin{proof}It is well known that $\beta_t$ is an isomorphism in this case. Also  
$\Gamma_g$ is in class ${\mathcal C}'$ and has property RD, so that  
$\mu_\sigma^\cA$ is an isomorphism. The result now follows from  
Proposition~\ref{K1twist}.
\end{proof}
\begin{remark}This corollary shows that we have obtained in particular  
the explicit generators for $K_1(C^*_r(\Gamma_g,\sigma))$, since $\Gamma_g$  
has property RD.
More explicitly, consider the standard presentation of $\Gamma_g$ in  
terms of generators and relations,
namely,
$$\Gamma_g = \left\{ a_j, b_j : \prod_{j=1}^g [a_j, b_j] = 1\right\}.$$
Then by Corollary \ref{thm:rs}, the unitary operators $\{T_{a_j}, T_{b_j}  \in U(\ell^2(\Gamma)):  
j=1,\ldots g\}$ form a natural
set of generators for  $K_1(C^*_r(\Gamma,\sigma))$ over $\mathbb Z$.
The corollary at the same time gives explicit generators for
 $K_1(\cA(\Gamma_g,\sigma))$ for any good unconditional completion.
\end{remark}

\subsection{$K$-homology}
Thanks to \cite{CHMM}, $K_1(C^*_r(\Gamma,\sigma)) \cong \mathbb Z^{2g}$
and $K_0(C^*_r(\Gamma,\sigma)) \cong \mathbb Z^{2}$  whenever $\Gamma =  
\Gamma_g$
as above. Moreover, we have determined a natural set of generators for   
$K_1(C^*_r(\Gamma_g,\sigma))$.
It was also shown in  \cite{CHMM} that $C^*_r(\Gamma,\sigma)$ is a  
$K$-amenable $C^*$-algebra. These
two facts together with the universal coefficient theorem \cite{RS}  
enable us to also compute the $K$-homology groups as
$$K^1(C^*_r(\Gamma,\sigma)) \cong \mathbb Z^{2g}, $$
and also
$$K^0(C^*_r(\Gamma,\sigma)) \cong \mathbb Z^{2},$$
where the $K$-homology groups $K^i(C^*_r(\Gamma,\sigma)) $ are defined  
as
usual as the Kasparov groups $KK^i(C^*_r(\Gamma,\sigma), \mathbb C)$  
for $i=0,1$.

%%%%%%%%%%%%%%%%%%%%%%%%%%%%%%%%%%%%%%%%%%%%%
\section*{Appendix: Twisted Rapid Decay}\label{RD}
\begin{center}by Indira Chatterji,\end{center} 
Department of Mathematics, Cornell University, Ithaca NY 14853, USA.\\
email: {indira@math.cornell.edu}
%%%%%%%%%%%%%%%%%%%%%%%%%%%%%%%%%%%%%%%%%%%%%%%
%\maketitle

Throughout this appendix, $\Ga$ is a finitely generated group, endowed with a length function $\ell$, and $\sigma$ is a multiplier on $\Ga$. We adopt the notations used in the first paragraph of the paper.
\begin{definition}\label{sRD} We will say that the group $\Gamma$ has  
\emph{$\sigma$-twisted Rapid Decay property (with respect to the length  
$\ell$)} if
$$H^{\infty}_{{\ell}}(\Gamma,\sigma)\subseteq C^*_r(\Gamma,\sigma).$$
We just say that the group $\Gamma$ has the \emph{Rapid Decay property  
(with respect to the length $\ell$)}, if it has the $\sigma$-twisted Rapid Decay 
property (with respect to the length $\ell$) for the constant multiplier  
1. For short, we shall say that a group $\Gamma$ has \emph{property $\sigma$-RD}  
if there esists a length function $\ell$ with respect to which $\Gamma$  
has the $\sigma$-twisted Rapid Decay property.
\end{definition}
\begin{remark}
In the context  of noncommutative geometry, 
the reduced $C^*$-algebra $C^*_r(\Gamma,\sigma)$ represents the space of {\em continuous} 
functions on a noncommutative manifold, and $H^{\infty}_{{\ell}}(\Gamma,\sigma)$
the space of of {\em smooth} 
functions on the same noncommutative manifold. 
This comes from the abelian case, where using Fourier transforms, one easily sees that $C^*_r(\bbZ^n)  \cong C(\mathbb T^n)$ and that
$H^{\infty}_{\ell}(\bbZ^n) \cong C^\infty(\mathbb T^n)$ (for the word length associated to the generating set $S=\{(\pm 1,0,\dots),\dots,(0,\dots,\pm 1)\}$ of $\bbZ^n$). The ($\sigma$-twisted) Rapid Decay property can be rephrased as the 
desirable property that every smooth function on the noncommutative
manifold is also a continuous function. 
\end{remark}
\begin{proposition}\label{prop:RDequiv}
Let $\sigma$ be a multiplier on $\Gamma$ and $\ell$  
be a length function on $\Gamma$. The following are equivalent:
\begin{itemize}
\item[(1)] $\Gamma$ has $\sigma$-twisted Rapid Decay (with respect to  
the length $\ell$).
\item[(2)] There exist constants $C,s>0$ such that for any  
$f\in\bbC(\Gamma,\sigma)$
$$\|f\|_{op}\leq C\|f\|_s.$$
\item[(3)]There exists a polynomial $P$ such that for any
$f\in\bbC(\Gamma,\sigma)$ and $f$ supported in a ball of radius $r$
$$\|f\|_{op}\leq P(r)\|f\|_{\ell^2\Gamma}.$$
\item[(4)] There exists a polynomial $P$ such that for any
$f,g\in\bbC(\Gamma,\sigma)$ and $f$ supported in a ball of radius $r$
$$\|f*_{\sigma}g\|_{\ell^2\Gamma}\leq P(r)\|f\|_{\ell^2\Gamma}\|g\|_{\ell^2\Gamma}.$$
\end{itemize}\end{proposition}
%%%
\begin{proof} $(1)\Leftrightarrow (2)$ As in the case of untwisted Rapid Decay, the inclusion $H^{\infty}_{{\ell}}(\Gamma,\sigma)\subseteq C^*_r(\Gamma,\sigma)$ is continuous since both inclusions $H^{\infty}_{{\ell}}(\Gamma,\sigma)\subseteq\ell^2\Gamma$ and $C^*_r(\Gamma,\sigma)\subseteq\ell^2\Gamma$ are continuous. Since $H^{\infty}_{{\ell}}(\Gamma,\sigma)$ is a Fr\'echet space, the continuity of the inclusion $H^{\infty}_{{\ell}}(\Gamma,\sigma)\subseteq C^*_r(\Gamma,\sigma)$ rephrases as the statement of (2). The converse is obvious since $H^{s+1}_{{\ell}}(\Ga)\subseteq H^s_{{\ell}}(\Ga)$.

$(2)\Rightarrow(3)\Rightarrow(4)$ Take $f\in\bbC(\Gamma,\sigma)$ supported in a ball  
of radius $r$, then
$$\|f\|_{op}\leq  
C\|f\|_s=C\sqrt{\sum_{\gamma\in\Gamma}|f(\gamma)|^2(1+\ell(\gamma))^{2s} 
}\leq C(1+r)^s\|f\|_{\ell^2\Gamma}.$$
Hence (3) follows. Since  
$\|f\|_{op}=\sup\{\frac{\|f*_{\sigma}g\|_{\ell^2\Gamma}}{\|g\|_{\ell^2\Gamma}}|{0\not=g\in\ell^2\Gamma}\}$ we deduce (4) as well. That (4) implies (3) is by definition of the operator norm.

$(3)\Rightarrow(2)$ For $n\in\bbN$, denote by  
$S_n=\{\gamma\in\Gamma|n\leq\ell(\gamma)<n+1\}$ the sphere of radius  
$n$. For $f\in\bbC(\Gamma,\sigma)$ we have:
$$\|f\|_{op}=\|\sum_{n=0}^{\infty}\lambda_{\sigma}(f|_{S_n})\|_{op}\leq\sum_{n=0}^{\infty}\|f|_{S_n}\|_{op},$$
so that using (3) we get the following bound
\begin{eqnarray*}\|f\|_{op}&\leq  
&\sum_{n=0}^{\infty}P(n+1)\|f|_{S_n}\|_{\ell^2\Gamma}\leq\sum_{n=0}^{\infty}C(n+1)^k\ 
|f|_{S_n}\|_{\ell^2\Gamma}\\
&\leq &  
C\sqrt{\sum_{n=0}^{\infty}(n+1)^{- 
2}}\sqrt{\sum_{n=0}^{\infty}(n+1)^{2k+2}\|f|_{S_n}\|_{\ell^2\Gamma}^2}\leq  
C'\|f\|_{k+1}\end{eqnarray*}
where $C'$ is some constant bigger than $C\pi/6$.
\end{proof}
The following proposition was known by Ji and Schweitzer \cite{tuile}, but the proof we give here might be shorter.
\begin{lemma}\label{1=>sigma}Let $\ell$ be a length function on
$\Gamma$. If $\Gamma$ has Rapid Decay (with respect to the length
$\ell$), then $\Gamma$ has $\sigma$-twisted Rapid Decay (with respect
to the length $\ell$) for any multiplier $\sigma$.\end{lemma}
%%%
\begin{proof}Take $\gamma\in\Gamma$, then:
$$|f*_{\sigma}g(\gamma)|=|\sum_{\mu\in\Gamma}f(\gamma^{- 
1}\mu)g(\mu)\sigma(\gamma^{- 
1}\mu,\mu)|\leq\sum_{\mu\in\Gamma}|f(\gamma^{- 
1}\mu)|\,|g(\mu)|=|f|*|g|(\gamma)$$
so that summing and squaring over $\gamma\in\Gamma$ yields
$$\|f*_{\sigma}g\|_{\ell^2\Gamma}\leq\||f|*|g|\|_{\ell^2\Gamma}\leq P(r)\|f\|_{\ell^2\Gamma}\|g\|_{\ell^2\Gamma}$$
and we conclude that $\Gamma$ has $\sigma$-twisted Rapid Decay using  
the previous proposition.
\end{proof}
The following corollary is the first part of Proposition 2.1 in \cite{LaffRD} with an obvious modification.
\begin{cor}[Noncommutative Sobolev Embedding Theorem]\label{remarque} 
Let $\ell$ be a length function on $\Gamma$.  
If $\Gamma$ has Rapid Decay (with respect to the length $\ell$), then  
there is a constant $S$ sufficiently large such that for any  multiplier $\sigma$ on  
$\Gamma$ and any $s\geq S$, $H^s_{\ell}(\Gamma,\sigma)$ is a Banach algebra
such that $H^s_{\ell}(\Gamma,\sigma) \subseteq C^*_r(\Gamma, \sigma)$.
\end{cor}
%%%
\begin{proof} Let $s$ be bigger than the degree of the polynomial of point (3) in Proposition \ref{prop:RDequiv}. We first have to show that there is a constant $K=K(s)$ such that for  
any $f,g\in\bbC(\Gamma,\sigma)$, $\|f*_{\sigma}g\|_s\leq  
K\|f\|_s\|g\|_s$. But this is true since  
$\|f*_{\sigma}g\|_s\leq\||f|*|g|\|_s$ and $\||f|*|g|\|_s\leq K'\|f\|_s\|g\|_s$ by Proposition 2.1 part (a) in \cite{LaffRD} (see also Proposition 8.15 in \cite{Valette}) since we assumed that $\Gamma$ has Rapid Decay (with respect to the length $\ell$). Therefore $H^s_{\ell}(\Gamma,\sigma)$ is a Banach algebra. By Lemma \ref{1=>sigma}, we know that since $\Gamma$ has property RD,  $\Gamma$ has property $\sigma$-twisted RD for any multiplier $\sigma$ on $\Gamma$, and hence $H^s_{\ell}(\Gamma,\sigma) \subseteq C^*_r(\Gamma, \sigma)$ follows from Proposition \ref{prop:RDequiv} part (2).
\end{proof}
%%%%%%
\begin{remark}
In the context of noncommutative geometry, 
Corollary \ref{remarque}  can be viewed as the 
analog of the {Sobolev Embedding Theorem}
for a compact manifold $M$, a simplified version of which saying that any function 
in the Sobolev space $W^{s,2}(M)$ for $s> \dim M/2$  is actually continuous. 
Indeed, using Fourier transforms, on can see that $W^{s,2}({\mathbb T}^n)\simeq H^{s}_{{\ell}}(\bbZ^n)$ for the word length associated to the generating set $S=\{(\pm 1,0,\dots),\dots,(0,\dots,\pm 1)\}$ of $\bbZ^n$, and that $C^*_r(\mathbb Z)\simeq C(\mathbb T^n)$.
\end{remark}
%%%%%%%%
\begin{example}Groups having Rapid Decay notably include: Polynomial  
growth groups (Jolissaint \cite{Jolissaint}), free groups (Haagerup  
\cite{Haagerup}) and more generally Gromov hyperbolic groups  
(Jolissaint-de la Harpe \cite{Harpe}), cocompact lattices in $SL_3(F)$  
where $F$ is the $p$-adic field $\bbQ_p$, $\bbR,\bbC,\bbH$ or  
$E_{6(-26)}$, as well as finite products of rank one Lie groups (see  
Rammagge-Robertson-Steger \cite{RRS}, Lafforgue \cite{LaffRD} and  
\cite{indira}) and all lattices in a rank one Lie group, see \cite{withKim}.\end{example}
\noindent
{\bf Question}: Is it possible to find a group $\Gamma$ which doesn't
have Rapid Decay, but which has $\sigma$-twisted Rapid Decay for some
multiplier $\sigma$ on $\Gamma$ (or does the converse of Lemma  
\ref{1=>sigma} hold)?

\medskip

The following is the second part of Proposition 1.2 of \cite{LaffRD} with a trivial change. But we still recall Lafforgue's proof below for the sake of completeness. 
\begin{proposition}\label{KtheoryIso}
Let $\ell$ be a length function on $\Gamma$. 
If $\Gamma$ has Rapid Decay (with respect to the length $\ell$), then  
 for any multiplier $\sigma$ on $\Gamma$ and for $s$ sufficiently large (and also for $s=\infty$), 
the inclusion  
$H^{s}_\ell (\Gamma,\sigma)\hookrightarrow C^*_r(\Gamma,\sigma)$ induces an  
isomorphism in $K$-theory.
\end{proposition}
\begin{proof}
The idea of the proof is as follows. By Corollary \ref{remarque}, there exists $S>0$ and finite 
such that for any $s\geq S$, $H^s_\ell(\Gamma,\sigma) \subseteq C^*_r(\Gamma, \sigma)$, and since $\bbC(\Gamma, \sigma)\subseteq H^{s}_\ell(\Gamma,\sigma)$, 
it follows that $H^{s}_\ell(\Gamma,\sigma)$ is a dense $*$-subalgebra of 
$C^*_r(\Gamma, \sigma)$. All we have to show is that the inclusion $H^s_\ell(\Gamma,\sigma) \subseteq C^*_r(\Gamma, \sigma)$ is spectral, it then follows (see e.g. Proposition 8.14 of \cite{Valette}) that the inclusion  $H^{s}_\ell(\Gamma,\sigma)\hookrightarrow C^*_r(\Gamma,\sigma)$ induces an isomorphism in $K$-theory.

Now, for two number $s,t$ such that $S<t<s$ the first step is to show that $H^s_\ell(\Gamma,\sigma)$ is stable by holomorphic functional calculus in $H^t_\ell(\Gamma,\sigma)$. To do so, and since $H^s_\ell(\Gamma,\sigma)$ is dense in $H^t_\ell(\Gamma,\sigma)$, it is enough (see Remark 8.13 in \cite{Valette}) to prove that the spectral radius $\rho_s(f)$ of $f\in H^s_\ell(\Gamma,\sigma)$ is the same as $\rho_t(f)$, the one of $f\in H^t_\ell(\Gamma,\sigma)$, namely that
\begin{equation}\label{rayonSpectral}\lim_{n\to\infty}\|f^{*_{\sigma}n}\|_s^{1/n}=\lim_{n\to\infty}\|f^{*_{\sigma}n}\|_t^{1/n},\end{equation}
where for $n\in\bbN$ we set $f^{*_{\sigma}n}=\underbrace{f*_{\sigma}f*_{\sigma}\dots *_{\sigma}f}_{n}$. Notice that since $t<s$, then $\|\ \|_t\leq\|\ \|_s$ and hence $\rho_t(f)\leq\rho_s(f)$, so we only need to prove the other inequality. For $\gamma\in\Gamma$, using the triangle inequality one sees that
\begin{eqnarray*}|f^{*_{\sigma}n}(\gamma)|&\leq&\sum_{\gamma_1,\dots,\gamma_{n-1}\in\Gamma}|f(\gamma\gamma_1^{-1})||f(\gamma_1\gamma_2^{-1})|\dots |f(\gamma_{n-2}\gamma_{n-1}^{-1})| |f(\gamma_{n-1})|\\
&=&\sum_{\gamma_1\dots\gamma_{n}=\gamma}|f(\gamma_1)|\dots |f(\gamma_n)|\end{eqnarray*}
Therefore, using that $(1+\ell(\gamma))^{s-t}\leq n^{s-t}\sum_{i=1}^n(1+\ell(\gamma_i))^{s-t}$ if $\gamma_1\dots\gamma_{n}=\gamma$ (which follows easily from Lemma 1.1.4 (3) in \cite{Jolissaint}) we deduce that
$$\|f^{*_{\sigma}n}\|_s=\|(1+\ell)^{s-t}f^{*_{\sigma}n}\|_t\leq n^{s-t+1}K^{n-1}\|f\|_s\|f\|_t^{n-1},$$
where $K=K(t)$ is the constant in the proof of Corollary \ref{remarque}. Taking the $n$-th root and the limit shows that $\lim_{n\to\infty}\|f^{*_{\sigma}n}\|_s^{1/n}\leq K\|f\|_t$. Replacing $f$ by $f^{*_{\sigma}m}$ in the previous inequlity, taking the $m$-th root and the limit shows $\rho_s(f)\leq\rho_t(f)$. We can now show that $H^s_\ell(\Gamma,\sigma) \subseteq C^*_r(\Gamma, \sigma)$ is spectral, namely that for $f\in H^s_\ell(\Gamma,\sigma)$, its spectral radius $\rho_s(f)$ equals $\rho_*(f)$, its spectral radius as an element of $C^*_r(\Gamma, \sigma)$. If $\rho_s(f)=0$ it is clear because $\rho_*(f)\leq\rho_s(f)$. Otherwise, H\"older's inequality shows that
$$\|f\|_t\leq\|f\|_s^{\frac{t}{s}}\|f\|_{\ell^2\Gamma}^{\frac{1-t}{s}},$$
and hence 
$$|f^{*_{\sigma}n}\|_{op}\geq\|f^{*_{\sigma}n}\|_{\ell^2\Gamma}\geq\|f^{*_{\sigma}n}\|^{\frac{s}{s-t}}_t\|f^{*_{\sigma}n}\|_s^{-\frac{t}{s-t}},$$
so that we conclude using equality (\ref{rayonSpectral}).\end{proof}
%%%%%%%%%%%%%%%%%

%%%%%%%%%%%%%%%%%%%%%%%%%%%


\begin{thebibliography}{CHMMM}
%%%
\bibitem[At]{At}
M.\ F.\ Atiyah.
\newblock {Elliptic operators, discrete groups and von Neumann  
algebras.}
\newblock Ast\'erisque {\bf 32/33} {(1976)}, 43--72.

%%

\bibitem[AtSi]{AtSi}
M.\ F.\ Atiyah, I. M. Singer.
\newblock \emph{The index of elliptic operators. III.}
\newblock Ann. of Math. (2) {\bf 87}, {(1968)}, 546--604.

%%

\bibitem[BaCo]{BC}
P. Baum, A. Connes.
\newblock \emph{$K$-theory for Lie groups and foliations.}
\newblock Enseign. Math. (2) {\bf 46} {(2000)}, no. 1-2, 3--42.

%%%
\bibitem[BaDo]{bd}
P. Baum and R. Douglas.
\newblock\emph{$K$-homology and index theory.}
\newblock Proceedings of Symposia in Pure Mathematics, {\bf 38}, Part 1  
(1982), 117--173.
%%%%

\bibitem[BeMaVa]{Bettaieb-Matthey} 
H.\ Bettaieb, M.\ Matthey, A.\  Valette.
\newblock\emph{Low-dimensional group homology and the Baum-Connes  
assembly map.}
\newblock Preprint \emph{1999}.

%%%%%%%%%%%

\bibitem[BrRo]{BrRo}  O. Bratteli,  D. Robinson.
\newblock Operator algebras and quantum statistical mechanics. 2.  
Equilibrium states. Models in quantum statistical mechanics.
\newblock Second edition. Texts and Monographs in Physics.  
Springer-Verlag, Berlin, 1997.

%%%

\bibitem[Bro]{Brown} 
K. Brown. 
\newblock \emph{Cohomology of groups.}
Graduate texts in Mathematics, vol. {\bf 87}, 
Springer-Verlag, New York, 1982.
%%%

\bibitem[BrSu]{BrSu}
J. Br\"uning and T. Sunada.
\newblock\emph{On the spectrum of gauge-periodic elliptic operators.}  
M\'ethodes semi-classiques, Vol. 2 (Nantes, 1991). {\em Ast\'erisque}  
{\bf 210} (1992), 65--74. {\em eidem.} On the spectrum of periodic  
elliptic operators,{\em Nagoya Math. J.} {\bf  126} (1992), 159--171.
%%%

\bibitem[CHMM]{CHMM}
A. Carey, K. Hannabuss, V. Mathai and P. McCann,
\newblock\emph{Quantum Hall Effect on the hyperbolic plane.}
\newblock Commun. Math. Phys. {\bf 190} No. 3 (1998) 629--673.
%%

\bibitem[Ch]{indira}
I. Chatterji.
\newblock \emph{Property (RD) for uniform lattices in products of rank  
one Lie groups with some rank two Lie groups.}
\newblock Geometria Dedicata {\bf 96} (2003) 161-177.
%%%

\bibitem[ChR]{withKim}
I. Chatterji, K. Ruane.
\newblock\emph{Some geometric groups with Rapid Decay}
\newblock {\tt [math.GR/0310356]}.

%%%

\bibitem[aT]{les_welches}
P.\ -A.\ Cherix, M.\ Cowling, P.\ Jolissaint, P.\ Julg, A.\ Valette.
\newblock \emph{Groups with the Haagerup property (Gromov's  
a-T-menability).}
\newblock Birk\"auser. Progress in Math. 197, 2001.
%%%

\bibitem[Co81]{Co81} A. Connes,
\newblock \emph{An analogue of the Thom
isomorphism for crossed products of a $C^*$ algebra by an action
of $\mathbb R$,}
  \newblock {Adv. in Math.} {\bf 39} (1981), 31-55.

%%%%%%%

\bibitem[CoMo]{ConnesMosc}
A. Connes, H. Moscovici.
\newblock \emph{Cyclic cohomology, the Novikov conjecture and  
hyperbolic groups.}
\newblock Topology {\bf 29} (1990), no. 3, 345--388.
%%%

\bibitem[G]{Greiner}
P. Greiner.
\newblock \emph{An asymptotic expansion for the heat equation.}
\newblock Arch. Ration. Mech. and Anal., {\bf 41} (1971), 163-218.
%%%

\bibitem[Gr]{Gr}
M. Gromov.
\newblock \emph{Volume and bounded cohomology.}
\newblock Publ. IHES., {\bf 56} (1983) 213-307.

%%%

\bibitem[Ha]{Haagerup}
U.\ Haagerup.
\newblock\emph{An example of a nonnuclear C*-algebra which has the  
metric approximation property.}
\newblock Invent. Math. {\bf 50} (1979) 279--293.
%%

\bibitem[dH]{Harpe}
P.\ de la Harpe.
\newblock \emph{Groupes hyperboliques, alg\`ebres d'op\'erateurs et un  
th\'eor\`eme de Jolissaint.}
\newblock C. R. Acad. Sci. Paris S\'er. I {\bf 307} {(1988)},  
771--774.
%%%

\bibitem[Hir]{Hir}
F. Hirzebruch.
\newblock \emph{Topological methods in algebraic geometry.}
\newblock Die Grundlehren der Mathematischen Wissenschaften, Band 131  
Springer-Verlag New York, Inc., New York 1966.

%%
\bibitem[Ji]{Ji}
R. Ji.
\newblock\emph{Smooth dense subalgebras of reduced group  
$C^*$-algebras, Schwartz cohomology of groups, and cyclic cohomology.}
\newblock J. Funct. Anal {\bf 107} (1992), 1--33.

%%
\bibitem[JiSc]{tuile}
R. Ji, L. B. Schweitzer.
\newblock \emph{Spectral invariance of smooth crossed products, and  
rapid decay locally compact groups.}
\newblock $K$-Theory {\bf 10} (1996), no. 3, 283--305.
%%

\bibitem[Jo]{Jolissaint}
P.\ Jolissaint.
\newblock \emph{Rapidly decreasing functions in reduced C*-algebras of  
groups.}
\newblock Trans. Amer. Math. Soc. {\bf 317} {(1990)}, 167--196.

%%
\bibitem[Ka]{Kasparov}
G.\ Kasparov.
\newblock \emph{$K$-theory, group C*-algebras, and higher signatures  
(Conspectus).}
\newblock Novikov conjectures, index theorems and rigidity, Vol. 1  
(Oberwolfach, 1993), 101--146,  London Math. Soc. Lecture Note Ser.,  
226,
\newblock Cambridge Univ. Press, 1995.

%%%%%%%%%%%%%%%
\bibitem[Ka2]{Ka2}
G.\ Kasparov.
\newblock \emph{Operator $K$-theory and its applications: elliptic operators, 
group representations, higher signatures, $C^*$-extensions. }
\newblock Proceedings of the International Congress of Mathematicians, 
Vol. 1, 2 (Warsaw, 1983), 987--1000, 
\newblock PWN, Warsaw, 1984.

%%%%%%%%%%%%%%%%%%%%%%%

\bibitem[Ko]{Ko} Yu. Kordyukov, $L^p$-theory of elliptic differential
operators on manifolds of bounded geometry, {\em Acta Appl.
Math.}, {\bf 23} (1991), {223--260}.


\bibitem[La]{Laff}
V.\ Lafforgue.
\newblock \emph{$K$-Th\'eorie bivariante pour les alg\`ebres de Banach  
et conjecture de Baum-Connes.}
\newblock Invent. Math. {\bf 149} (2002), no. 1, 1--95.
%%%

\bibitem[La2]{LaffRD}
V.\ Lafforgue.
\newblock \emph{A proof of property (RD) for discrete cocompact  
subgroups of $SL_3({\bf R})$ and $SL_3({\bf C})$.}
\newblock Journal of Lie Theory {\bf 10}, {(2000)}, 255--267.
%%

\bibitem[Ma]{Ma}
V. Mathai.
\newblock \emph{$K$-theory of twisted group $C^*$-algebras and positive  
scalar curvature.}
\newblock Contemp. Math. {\bf 231} (1999), 203--225.
%%

\bibitem[Ma2]{Ma2}
V. Mathai.
\newblock \emph{On positivity of the Kadison constant and  
noncommutative Bloch theory.}
\newblock Tohoku Mathematical Publications, {\bf 20} (2001) 107-124.
%%

\bibitem[MaMa]{MM}
M. Marcolli, V. Mathai.
\newblock \emph{Twisted index theory on good orbifolds, I:  
noncommutative Bloch theory.}
\newblock Communications in Contemporary Mathematics, {\bf 1} (1999)  
553-587.

%%
\bibitem[MiYu]{Mineyev_Yu}
I. Mineyev, G. Yu.
\newblock \emph{The Baum-Connes conjecture for hyperbolic groups.}
\newblock Invent. Math. {\bf 149} (2002) 1-95.
%%

\bibitem[Na]{Natsume}
T.\ Natsume.
\newblock \emph{The Baum-Connes conjecture, the commutator theorem and  
Rief\-fel projections.}
\newblock C. R. Math. Rep. Acad. Sci. Canad. {\bf 1} {(1988)},  
13--18.

%%

\bibitem[PR]{PR} J. Packer, I. Raeburn, Twisted cross products of
$C^*$-algebras,
{\em Math.\ Proc.\ Camb.\ Phil.\ Soc.} {\bf 106} (1989), 293-311.

%%

\bibitem[PiVo]{PV}
M. Pimsner and D. Voiculescu.
\newblock \emph{Exact sequences for $K$-groups and Ext-groups of  
certain cross-product $C\sp{*}$-algebras.}
\newblock J. Operator Theory, {\bf 4} (1980), no. 1, 93-118.

%%
\bibitem[RaRoSt]{RRS}
J.\ Ramagge, G.\ Robertson, T.\ Steger.
\newblock \emph{A Haagerup inequality for $\widetilde A\sb  
1\times\widetilde A\sb 1$ and $\widetilde A\sb 2$ buildings.}
\newblock Geom. Funct. Anal. {\bf 8 }{(1998)}, no. {\bf 4}, 702--731.

%%
\bibitem[Ri]{Rieff}
M. Rieffel.
\newblock \emph{$C^*$-algebras associated with irrational rotations.}
\newblock Pac. J. Math. {\bf 93} (1981), 415-429.

%%
\bibitem[RS]{RS}
J. Rosenberg, C. Schochet,
\newblock \emph{The K\" unneth theorem and the universal coefficient  
theorem for Kasparov's generalized $K$-functor.}
\newblock Duke Math. J. {\bf 55} (1987), no. 2, 431-474.

%%
\bibitem[Va]{Valette} A. Valette.
\newblock\emph{Introduction to the Baum-Connes Conjecture.} Notes taken  
by Indira Chatterji. With an appendix by Guido Mislin. Lectures in  
Mathematics ETH Z\"urich.
Birkh\"auser Verlag, Basel, 2002.


\end{thebibliography}
\end{document}